\documentclass[12pt]{article}
\RequirePackage{amsthm,amsmath,natbib}
\usepackage{amsmath, amsthm, amssymb}
\usepackage{verbatim,color,amssymb,epsfig}
\usepackage[usenames,dvipsnames]{xcolor}
\usepackage{multirow,booktabs}

\usepackage{xr}
\RequirePackage{amsthm,amsmath,natbib}
\usepackage{amsmath, amsthm, amssymb}
\usepackage[auth-sc, affil-sl]{authblk}
\usepackage{epsfig,natbib}

\renewcommand{\baselinestretch}{1.3}
\renewcommand{\baselinestretch}{1.5}
\def\boxit#1{\vbox{\hrule\hbox{\vrule\kern6pt
          \vbox{\kern6pt#1\kern6pt}\kern6pt\vrule}\hrule}}

\def\trans{^{\rm T}}

\numberwithin{equation}{section}
\theoremstyle{plain}

\newtheorem{Theorem}{\underline{\bf Theorem}}[section]

\newtheorem{Remark}{\underline{\bf Remark}}

\newtheorem{Condition}{\underline{\bf Condition}}
\newtheorem{Definition}{\underline{\bf Definition}}

\begin{document}

\thispagestyle{empty}
\baselineskip=28pt
\begin{center}
{\LARGE{\bf Asymptotic optimality of sparse linear discriminant analysis with arbitrary number of classes}}
\end{center}
\vskip 10mm

\baselineskip=14pt
\vskip 2mm
\begin{center}
 Ruiyan Luo\\
\vskip 2mm
Division of Epidemiology and Biostatistics, Georgia State University School of Public Health, One Park Place, Atlanta, GA 30303\\
rluo@gsu.edu\\
\hskip 5mm\\
Xin Qi\\
\vskip 2mm
Department of Mathematics and Statistics, Georgia State University, 30 Pryor Street, Atlanta, GA 30303-3083\\
xqi3@gsu.edu \\
\end{center}

\begin{abstract}
Many sparse linear discriminant analysis (LDA) methods have been proposed to overcome the major problems of the classic LDA in high-dimensional settings. However, the asymptotic optimality
results are limited to the case that there are only two classes, which is due to the fact that the classification boundary of LDA is a hyperplane and  explicit formulas exist for the classification error in this case. In the situation where there are more than two classes, the classification boundary is usually complicated and no explicit formulas for the classification errors exist. In this paper, we consider the asymptotic optimality in the high-dimensional settings for a large family of linear classification rules with arbitrary number of classes   under the situation of multivariate normal distribution. Our main theorem provides easy-to-check criteria for the asymptotic optimality of a general classification rule in this family as dimensionality and sample size both go to infinity and the number of classes is arbitrary. We   establish the corresponding convergence rates. The general theory is applied to the classic LDA and the extensions of two recently proposed sparse LDA methods to obtain the asymptotic optimality.  We conduct simulation studies on the extended methods in various settings.
\end{abstract}

 \vspace{6mm}
\par\medskip\noindent
\underline{\bf Short title}: Asymptotic optimality of sparse LDA

\pagenumbering{arabic}
\newlength{\gnat}
\setlength{\gnat}{22pt}
\baselineskip=\gnat

\section{Introduction}\label{sec1}
As an important classification method, the linear discriminant analysis (LDA) performs well in the settings of small $p$ and large $n$. However, it faces major problems for high-dimensional data with large $p$ and small $n$. In theory, \citet{bickel2004some} and \citet{Shao-2011} showed that in the case of $p>n$, the classic LDA can be as bad as the random guessing. To address these problems, various regularized discriminant analysis methods have been proposed, including those described in \citet{Friedman-1989}, \citet{Krzanowski-1995}, \citet{Dudoit-2001}, \citet{bickel2004some}, \citet{Guo-2007}, \citet{Xu-2009}, \citet{Tibshirani-2009}, \citet{Witten-2011}, \citet{Clemmensen-2011}, \citet{Shao-2011}, \citet{cai2011direct}, \citet{fan2012road},  and many others. Asymptotic optimality has been established in some of these papers \citep{Shao-2011, cai2011direct, fan2012road}. However, these asymptotic optimality
results are limited to the case where there are only two classes. 
The major difficulty preventing the derivation of asymptotic optimality for the multiclass classification is that for the two-class classification, the classification boundary of LDA is a hyperplane and there exist explicit formulas for the classification error, however, when the number of classes is greater than two, the classification boundary is usually complicated and no explicit formula for the classification error exists.

In this paper, we consider the asymptotic optimality in high-dimensional settings for a large family of linear classification rules with arbitrary number of classes   under the situation of multivariate normal distribution. The classification rules of the optimal LDA, the classic LDA, and those in \citet{Shao-2011} and \citet{cai2011direct} all belong to this family.  We first provide an upper bound on the difference between the conditional classification error of any classification rule in this family and the optimal classification error  for arbitrary  $n$,  $p$ and $K$ (the number of classes). Through an example, we illustrate that there exist situations where this bound is asymptotic optimal. Based on the upper bound, we develop our main theorem which provides the conditions leading to the asymptotic optimality for a general classification rule in this family as dimensionality and sample size both go to infinity and the number of classes is arbitrary. These conditions are relatively easily verified for various particular classification rules in this family. We establish the convergence rates for the asymptotic optimality under these conditions. Then we apply this  theorem to several  particular classification rules. 
   We extend the sparse LDA methods in \citet{Shao-2011} and \citet{cai2011direct} from the two-class situations to the multi-class situations, and apply our general theorem to the two extended methods to obtain the asymptotic optimality and the corresponding convergence rates. Simulation studies are performed to evaluate the predictive performance of the two extended methods in various settings.

The rest of this paper is organized as follows. In Section 2, after introducing the notations and main assumptions, the classic and two sparse LDA methods are described shortly. In Section 3, we introduce a  family of linear classification rules and provide the main theorems.  A necessary condition for the asymptotic optimality of the usual LDA and corresponding convergence rate are given in Section 4 as $p$, $n$ and $K$ all go to infinity. In Sections 5 and 6, two sparse LDA methods in \citet{Shao-2011} and \citet{cai2011direct} are extended to the multiclass cases, and the asymptotic optimality and the corresponding convergence rates are provided. Simulation studies are performed in Section 7. A short discussion is provided in Section 8. All the proofs can be found in the supplementary materials.

\section{LDA and sparse LDA}\label{sec2}
We first introduce some notations. For any vector $\mathbf{v}=(v_1,\cdots,v_p)\trans$, let $\|\mathbf{v}\|_2$, $\|\mathbf{v}\|_1$ and $\|\mathbf{v}\|_{\infty}=\max_{1\le i\le p}|v_i|$ be the $l_2$, $l_1$ and $l_{\infty}$ norms of $\mathbf{v}$, respectively. For any $p\times p$ symmetric nonnegative definite matrix $\mathbf{M}$, we use $\lambda_{max}(\mathbf{M})$  and $\lambda_{min}(\mathbf{M})$ to denote its largest  and smallest eigenvalues, respectively.  We define two norms for $\mathbf{M}$, 
\begin{align*}
\|\mathbf{M}\|=\sup_{\mathbf{v}\in\mathbb{R}^p, \|\mathbf{v}\|_2=1}\|\mathbf{M}\mathbf{v}\|_2 ,\quad  \text{ and }\quad  \|\mathbf{M}\|_{\infty}=\max_{1\le k,l\le p}|M_{kl}|
\end{align*}
where $M_{kl}$ is the $(k,l)$th entry of $\mathbf{M}$.  The first one is  the {\it operator norm} and the second one is the {\it max norm}.

In this paper, we assume that there are $K$ classes and the population in the $i$th class has a multivariate normal distribution $N_p(\boldsymbol{\mu}_i,  \boldsymbol{\Sigma})$, where $\boldsymbol{\mu}_i$ is the $i$th class mean, $1\le i\le K$, and $\boldsymbol{\Sigma}$ is the common covariance matrix for all classes. We assume that the prior probabilities for all the classes are the same and equal to $1/K$. 
We will consider the situations where both $n$ and $p$ go to infinity, and  $K$ is arbitrary. 

We first present two key regularity conditions for our theory. 
\begin{Condition}\label{condition_1}
There is a constant $c_0$ (independent of p and K) such that
\begin{align*}
\frac{1}{c_0}\le  \lambda_{min}(\boldsymbol{\Sigma})\le\lambda_{max}(\boldsymbol{\Sigma})\le c_0. 
\end{align*}
\end{Condition}

\begin{Condition}\label{condition_2}
There exists a constant $c_1>0$ which does not depend on $p$ and $K$, such that
\begin{align*}
 \min_{1\le i\neq j\le K}(\boldsymbol{\mu}_i-\boldsymbol{\mu}_j)\trans\boldsymbol{\Sigma}^{-1}(\boldsymbol{\mu}_i-\boldsymbol{\mu}_j) \ge c_1.
\end{align*}
\end{Condition}
  Condition \ref{condition_1} is   the same as (2) in \citet{Shao-2011}.  To understand the meaning of the inequality in Condition  \ref{condition_2}, we consider the optimal linear discriminant rule.  For any $\mathbf{x}$, $\mathbf{y}\in \mathbb{R}^p$, let $d_{\boldsymbol{\Sigma}^{-1}}(\mathbf{x},\mathbf{y})=\sqrt{(\mathbf{x}-\mathbf{y})\trans\boldsymbol{\Sigma}^{-1}(\mathbf{x}-\mathbf{y})}$ which is the well known  Mahalanobis distance.  The optimal linear classification rule is given by  (see Section 6.8 in \citet{anderson2003introduction} or Theorem 13.2 in \citet{hardle2007applied}), 
\begin{align}
& T_{OPT} \text{: to allocate a new observation $\mathbf{x}$ to the $i$th class if  $d_{\boldsymbol{\Sigma}^{-1}}(\mathbf{x}, \boldsymbol{\mu}_i)$ is the smallest }\notag\\
 &\qquad\quad \text{among the $K$ distances: }  d_{\boldsymbol{\Sigma}^{-1}}(\mathbf{x}, \boldsymbol{\mu}_1), d_{\boldsymbol{\Sigma}^{-1}}(\mathbf{x}, \boldsymbol{\mu}_2),\cdots, d_{\boldsymbol{\Sigma}^{-1}}(\mathbf{x}, \boldsymbol{\mu}_K).\label{8}
\end{align}  
  We use $R_{OPT}$ to denote the misclassification rate of $T_{OPT}$. It is well known that under the assumptions on the population distributions in this paper, $T_{OPT}$ is the Bayes rule and $R_{OPT}$ is the smallest among the misclassification rates of all possible classification rules.  Condition  \ref{condition_2} implies that  the squared Mahalanobis distance between any two class means is not less than $c_1$. If this condition is not satisfied, some class means will approach each other as $n$, $p\to\infty$ and these classes will be completely mixed together. In the case of two classes, we have an explicit formula for $R_{OPT}$: $R_{OPT}=\Phi(-(\boldsymbol{\mu}_1-\boldsymbol{\mu}_2)\trans\boldsymbol{\Sigma}^{-1}(\boldsymbol{\mu}_1-\boldsymbol{\mu}_2)/2)$, where $\Phi$ is the cumulative distribution function of the standard normal distribution (see Section  13.1 in \citet{hardle2007applied}). If Condition  \ref{condition_2} is not satisfied, we have $R_{OPT}\to 1/2$, which is the misclassification rate of a random guess. Condition \ref{condition_2} excludes these situations.  
	
 In practice, $\boldsymbol{\mu}_i$, $1\le i\le K$, and $\boldsymbol{\Sigma}$ are all unknown. Let $\mathbf{X}=\{\mathbf{x}_{ij}: 1\le i\le K, 1\le j\le n_i\}$ be a sample data set from the population,  where $\mathbf{x}_{ij}$ is the $j$th observation from the $i$th class and $n_i$ is the number of the observations from the $i$th class. Throughout this paper, we use 
\begin{align*}
&\bar{\mathbf{x}}_i=\frac{1}{n_i}\sum_{j=1}^{n_i}\mathbf{x}_{ij}, \quad\widehat{\boldsymbol{\Sigma}}=\frac{1}{n}\sum_{i=1}^K\sum_{j=1}^{n_i}(\mathbf{x}_{ij}-\bar{\mathbf{x}}_i)(\mathbf{x}_{ij}-\bar{\mathbf{x}}_i)\trans, \quad 1\le i\le K,
\end{align*}
  to denote the sample class means and the sample within-class covariance matrix which are estimates of $\boldsymbol{\mu}_i$ and $\boldsymbol{\Sigma}$, respectively.
The classic LDA rule is given by
\begin{align*}
& T_{LDA} \text{: to allocate a new observation $\mathbf{x}$ to the $i$th class if $d_{\widehat{\boldsymbol{\Sigma}}^{-1}}(\mathbf{x}, \bar{\mathbf{x}}_i)$ is the smallest among}\notag\\
 &\qquad \qquad d_{\widehat{\boldsymbol{\Sigma}}^{-1}}(\mathbf{x}, \bar{\mathbf{x}}_1),d_{\widehat{\boldsymbol{\Sigma}}^{-1}}(\mathbf{x}, \bar{\mathbf{x}}_2),\cdots, d_{\widehat{\boldsymbol{\Sigma}}^{-1}}(\mathbf{x}, \bar{\mathbf{x}}_K). 
\end{align*}
Unlike $T_{OPT}$,  the rule $T_{LDA}$ depends on the sample $\mathbf{X}$. It has been argued in Chapter 1 in \citet{devroye2013probabilistic} that  the conditional misclassification rate is a more natural measure of the predictive performance of a classification rule built based on $\mathbf{X}$ than the unconditional misclassification rate. Let $T$ be any linear classification rule  based on $\mathbf{X}$. The conditional misclassification rate of $T$ given $\mathbf{X}$ is defined as
\begin{align*}
&R_{T}(\mathbf{X})= \sum_{i=1}^KP \left(\left\{\text{$\mathbf{x}_{\rm new}$ belongs to the $i$-th class but $T(\mathbf{x}_{\rm new})\neq i$}\right\}\bigg\vert\mathbf{X} \right)
\end{align*}  
where $\mathbf{x}_{\rm new}$ is a new observation independent of $\mathbf{X}$.   $R_{T}(\mathbf{X})$ is  a function of $\mathbf{X}$ and the unconditional misclassification rate is the expectation of $R_{T}(\mathbf{X})$ with respect to the distribution of $\mathbf{X}$.   For simplicity, we use $R_{{LDA}}(\mathbf{X})$ to denote the conditional misclassification rate of $T_{LDA}$.

In the high-dimensional settings, the classic LDA performs poorly and can even fail completely \citep{bickel2004some, Shao-2011}. To revise the classic LDA in the high-dimensional settings, we note that 
 $ d_{{\boldsymbol{\Sigma}}^{-1}}(\mathbf{x}, \boldsymbol{\mu}_i)<d_{{\boldsymbol{\Sigma}}^{-1}}(\mathbf{x}, \boldsymbol{\mu}_j)$  is equivalent to
\begin{align*}
(\boldsymbol{\mu}_j-\boldsymbol{\mu}_i)\trans\boldsymbol{\Sigma}^{-1}(\mathbf{x}-\frac{\boldsymbol{\mu}_i+\boldsymbol{\mu}_j}{2})<0.
\end{align*}  
\citet{Shao-2011} imposed sparsity conditions on ${\boldsymbol{\Sigma}}$ and $\boldsymbol{\delta}=\boldsymbol{\mu}_2-\boldsymbol{\mu}_1$, and proposed the following sparse LDA classification rule for two classes:
\begin{align}
T_{SLDA} &\text{: to allocate a new observation $\mathbf{x}$ to the first class if }  \widetilde{\boldsymbol{\delta}}\trans\widetilde{\boldsymbol{\Sigma}}^{-1}(\mathbf{x}- \frac{\bar{\mathbf{x}}_1+\bar{\mathbf{x}}_2}{2})<0,\notag\\
&\quad\text{ and to the second class otherwise, }\label{20}
\end{align}
where   $\widetilde{\boldsymbol{\Sigma}}$ and $\widetilde{\boldsymbol{\delta}}$ are the thresholding estimates of ${\boldsymbol{\Sigma}}$ and $\boldsymbol{\delta}$, respectively.   We use $R_{\text{SLDA}}(\mathbf{X})$ to denote the conditional misclassification rate of $T_{SLDA}$ given the sample $\mathbf{X}$.

\citet{cai2011direct}  further observed that in the case of two classes, the optimal classification rule $T_{OPT}$ depends on ${\boldsymbol{\Sigma}}$ only through $\boldsymbol{\beta}={\boldsymbol{\Sigma}}^{-1}\boldsymbol{\delta}$. Hence, they only assumed that $\boldsymbol{\beta}$ is sparse and proposed a sparse estimate $\widehat{\boldsymbol{\beta}}$ of  $\boldsymbol{\beta}$ based on a linear programming optimization problem.  Then they proposed the following linear programming discriminant (LPD) rule  for two classes,
\begin{align}
T_{LPD} &\text{: to allocate a new observation $\mathbf{x}$ to the first class if }  \widehat{\boldsymbol{\beta}}\trans (\mathbf{x}- \frac{\bar{\mathbf{x}}_1+\bar{\mathbf{x}}_2}{2})<0,\notag\\
&\quad\text{ and to the second class otherwise. }\label{21}
\end{align}
We use $R_{LPD}(\mathbf{X})$ to denote the conditional misclassification rate of $T_{LPD}$ given the  sample $\mathbf{X}$. The following definition of asymptotic optimality has been used in  \citet{Shao-2011}, \citet{cai2011direct} and other papers.
 
 \begin{Definition}\label{definition_1}
\renewcommand{\baselinestretch}{1.5}
 Let $T$ be a linear classification rule with conditional misclassification rate $R_T(\mathbf{X})$. Then $T$ is asymptotically optimal if 
\begin{align}
&\frac{R_{T}(\mathbf{X})}{R_{{OPT}}}-1= o_p(1).\label{112}
\end{align}
\end{Definition}

Since $0\le R_{{OPT}}\le R_{T}(\mathbf{X})\le 1$ for any $\mathbf{X}$, $\eqref{112}$ implies that $ 0\le R_{T}(\mathbf{X})-R_{{OPT}}= o_p(1)$. Hence we have $R_{T}(\mathbf{X})\to R_{{OPT}}$ in probability and $E[R_{T}(\mathbf{X})]\to R_{{OPT}}$, which have been used to define the consistency of a classification rule by  \citet{devroye2013probabilistic} and others. If $R_{{OPT}}$ is bounded away from 0, then $  R_{T}(\mathbf{X})-R_{{OPT}}= o_p(1)$ also implies $\eqref{112}$. However, if $R_{{OPT}}\to 0$, $\eqref{112}$ is stronger than  $  R_{T}(\mathbf{X})-R_{{OPT}}=o_p(1)$.

\section{Upper bounds and convergence rates}
In this section, we consider a family of linear classification rules motivated by the following observation. The optimal classification rule $T_{OPT}$ can be rewritten in the following way. Let
\begin{align}
\mathbf{a}_{ji}=\boldsymbol{\Sigma}^{-1/2}(\boldsymbol{\mu}_j-\boldsymbol{\mu}_i), \quad \mathbf{b}_{ji}=\frac{1}{2}(\boldsymbol{\mu}_j+\boldsymbol{\mu}_i),\label{22}
\end{align}  
where $1\le i, j\le K$. Then $T_{OPT}$ assigns a new observation $\mathbf{x}$ to the $i$th class if $\mathbf{a}_{ji}\trans\boldsymbol{\Sigma}^{-1/2}(\mathbf{x}-\mathbf{b}_{ji})<0$ for all $j\neq i$. Based on this observation, we consider the family of classification rules having the form, 
\begin{align}
T \text{: to assign a new  $\mathbf{x}$ to the $i$th class if } \widehat{\mathbf{a}}_{ji}\trans\boldsymbol{\Sigma}^{-1/2}(\mathbf{x}- \widehat{\mathbf{b}}_{ji})<0, \text{ for all } j\neq i, \label{23} 
\end{align}  
where $\widehat{\mathbf{a}}_{ji}$ and $\widehat{\mathbf{b}}_{ji}$ are  $p$-dimensional vectors which may depend on the sample $\mathbf{X}$, and  satisfy
\begin{align}
\widehat{\mathbf{a}}_{ji}=-\widehat{\mathbf{a}}_{ij},\quad \widehat{\mathbf{b}}_{ji}=\widehat{\mathbf{b}}_{ij},\label{24}
\end{align}
for all $1\le i\neq j\le K$. In addition to $T_{OPT}$, the classification rules $T_{LDA}$, $T_{SLDA}$ and $T_{LPD}$ all belong to this family. The specific expressions of $\widehat{\mathbf{a}}_{ji}$ and $\widehat{\mathbf{b}}_{ji}$ for them will be given in the following sections. We first provide an upper bound on $R_T(\mathbf{X})-R_{OPT}$ in the following theorem.   

\begin{Theorem}\label{theorem_1}
  
\begin{align}
&R_{T}(\mathbf{X})-R_{{OPT}}\label{25}\\
\le& \frac{1}{K}\sum_{i=1}^K\sum_{j\neq i}P\left(\widehat{\mathbf{a}}_{ji}\trans\mathbf{Z}>\widehat{\mathbf{a}}_{ji}\trans\boldsymbol{\Sigma}^{-1/2}(\widehat{\mathbf{b}}_{ji}-\boldsymbol{\mu}_i),  \quad\mathbf{a}_{ji}\trans\mathbf{Z}<\mathbf{a}_{ji}\trans\boldsymbol{\Sigma}^{-1/2}(\mathbf{b}_{ji}-\boldsymbol{\mu}_i) \bigg\vert\mathbf{X} \right),\notag
\end{align}
where $\mathbf{Z}$ is a $p$-dimensional random vector with distribution $N(\mathbf{0},\mathbf{I}_p)$ and independent of the  sample $\mathbf{X}$. For the special case of $K=2$, we have the following equality,
\begin{align}
&R_{T}(\mathbf{X})-R_{{OPT}}\label{100}\\
=& \frac{1}{2}\sum_{i=1}^2\sum_{j\neq i}\left[P \left(\widehat{\mathbf{a}}_{ji}\trans\mathbf{Z}>\widehat{\mathbf{a}}_{ji}\trans\boldsymbol{\Sigma}^{-1/2}(\widehat{\mathbf{b}}_{ji}-\boldsymbol{\mu}_i),  \quad\mathbf{a}_{ji}\trans\mathbf{Z}<\mathbf{a}_{ji}\trans\boldsymbol{\Sigma}^{-1/2}(\mathbf{b}_{ji}-\boldsymbol{\mu}_i)\bigg\vert\mathbf{X} \right) \right.\notag\\
&\qquad\left. -P \left(\widehat{\mathbf{a}}_{ji}\trans\mathbf{Z}<\widehat{\mathbf{a}}_{ji}\trans\boldsymbol{\Sigma}^{-1/2}(\widehat{\mathbf{b}}_{ji}-\boldsymbol{\mu}_i),  \quad\mathbf{a}_{ji}\trans\mathbf{Z}>\mathbf{a}_{ji}\trans\boldsymbol{\Sigma}^{-1/2}(\mathbf{b}_{ji}-\boldsymbol{\mu}_i)\bigg\vert\mathbf{X} \right) \right]\;.\notag
\end{align}
\end{Theorem}
\begin{Remark}\label{remark_4}
 In $\eqref{25}$, $P \left(\widehat{\mathbf{a}}_{ji}\trans\mathbf{Z}>\widehat{\mathbf{a}}_{ji}\trans\boldsymbol{\Sigma}^{-1/2}(\widehat{\mathbf{b}}_{ji}-\boldsymbol{\mu}_i),  \quad\mathbf{a}_{ji}\trans\mathbf{Z}<\mathbf{a}_{ji}\trans\boldsymbol{\Sigma}^{-1/2}(\mathbf{b}_{ji}-\boldsymbol{\mu}_i)\bigg\vert\mathbf{X} \right)$ denotes the conditional probability of the event $\{\widehat{\mathbf{a}}_{ji}\trans\mathbf{Z}>\widehat{\mathbf{a}}_{ji}\trans\boldsymbol{\Sigma}^{-1/2}(\widehat{\mathbf{b}}_{ji}-\boldsymbol{\mu}_i)\}\bigcap\{\mathbf{a}_{ji}\trans\mathbf{Z}<\mathbf{a}_{ji}\trans\boldsymbol{\Sigma}^{-1/2}(\mathbf{b}_{ji}-\boldsymbol{\mu}_i)\}$ given $\mathbf{X}$. Since $\mathbf{Z}$ is independent of $\mathbf{X}$, when we calculate the above conditional probability, we just need to consider $\widehat{\mathbf{a}}_{ji}$ and $\widehat{\mathbf{b}}_{ji}$ as constants and calculate the probability with respect to the distribution of $\mathbf{Z}$. The same interpretations will be applied  throughout this paper.  
\end{Remark}
In the case of $K=2$, we obtained an equality $\eqref{100}$, but for $K>2$, the upper bound on the right hand side of $\eqref{25}$ is usually greater than $R_{T}(\mathbf{X})-R_{{OPT}}$. So a natural question is whether the upper bound  is asymptotically optimal, that is, whether the ratio between $R_{T}(\mathbf{X})-R_{{OPT}}$ and the bound converges to 1. However, it is hard to answer this question for general situations. We provide an example where the ratio converges to 1. This example is somewhat artificial, but serves to illustrate that there exist situations where the bound in $\eqref{102}$ is asymptotically optimal and  gives some insights into the probabilities in the upper bound.
	\begin{figure}[h]
	\begin{center}
\includegraphics[height=4in,width=4in]{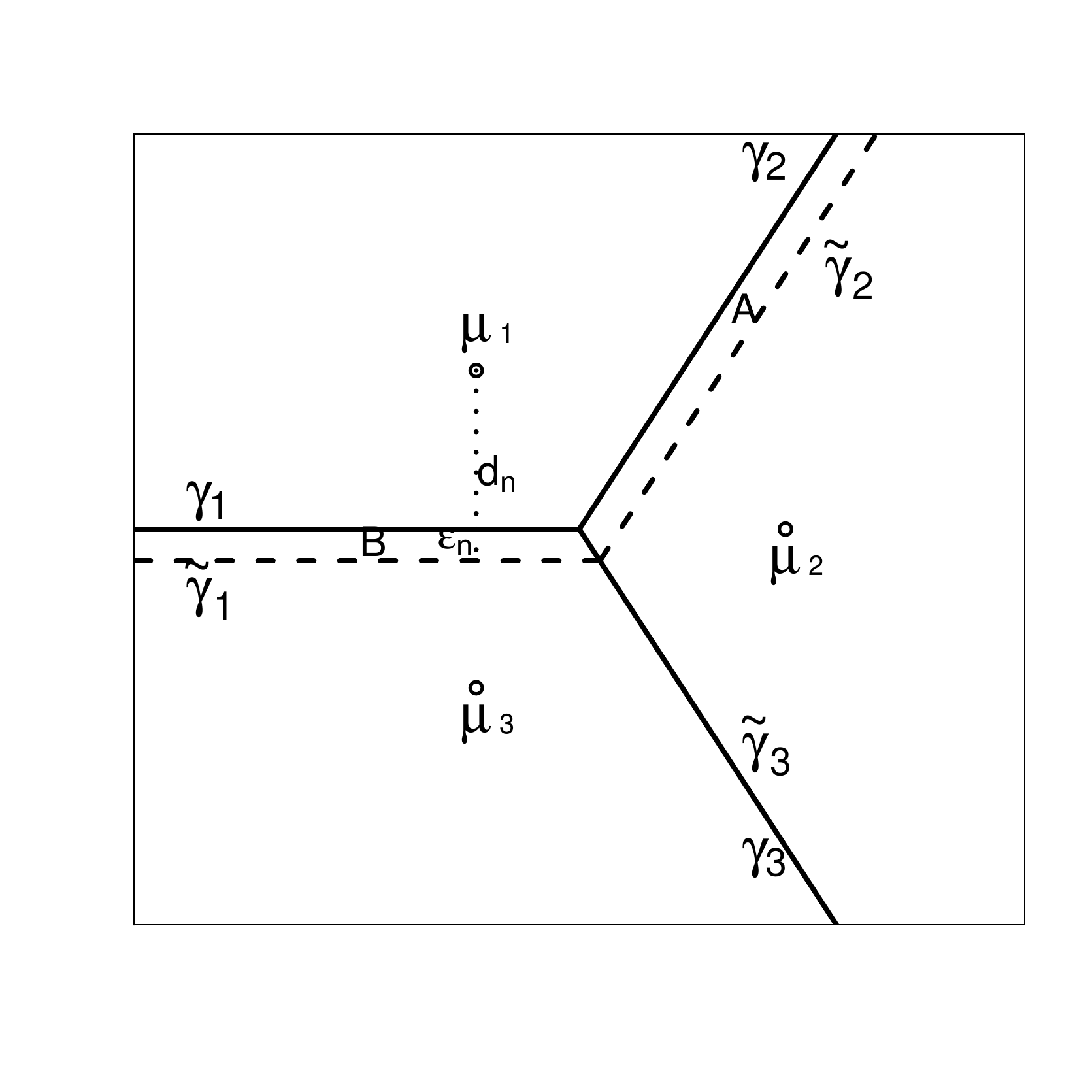}
    \caption{ \baselineskip=10pt Illustration of an example where the upper bound in $\eqref{25}$ is asymptotically optimal. Here K=3 and $\boldsymbol{\mu}_1$, $\boldsymbol{\mu}_2$ and $\boldsymbol{\mu}_3$ are the class means and the distance between any two of them is the same and equal to $2d_n$. The solid lines, $\gamma_1$, $\gamma_2$ and $\gamma_3$, denote the boundary of $T_{OPT}$ and the dashed lines, $\widetilde{\gamma_1}$, $\widetilde{\gamma_2}$ and $\widetilde{\gamma_3}$, are the classification boundary of $T$.  $\widetilde{\gamma_3}$ is superposed on $\gamma_3$. $d_n$ and $d_n+\epsilon_n$ are distances from $\boldsymbol{\mu}_1$ to the boundary lines of $T_{OPT}$ and $T$, respectively.  A is the region bounded by $\widetilde{\gamma_2}$, $\gamma_2$ and $\gamma_3$, and  B is the region bounded by $\widetilde{\gamma_1}$, $\gamma_1$ and $\gamma_3$, respectively. }
\label{fig_2}
\end{center}
\end{figure}
 The example is shown in Figure \ref{fig_2}, where $K=3$ and $p=2$ are fixed. We assume that the pairwise distance between any two pairs of class means are the same, and $\boldsymbol{\Sigma}=\bigl(\begin{smallmatrix} 1&0 \\ 0&1 \end{smallmatrix}\bigr)$. The boundary lines $\widetilde{\gamma_1}$ and $\widetilde{\gamma_2}$ are parallel to $\gamma_1$ and $\gamma_2$, respectively, and the distances between boundary lines are the same and equal to $\epsilon_n$. The regions A is  actually the collection of points which are assigned to class 2 under $T_{OPT}$ but assigned to class 1 under $T$, and B is the collection of points which are assigned to class 3 under $T_{OPT}$ but assigned to class 1 under $T$. Hence, we can see that 
\begin{align}
&R_{T}(\mathbf{X})-R_{{OPT}}=\frac{1}{3}[P_2(A)+P_3(B)-P_1(A)-P_1(B)]=\frac{2}{3}[P_3(B)-P_1(B)],\label{162}
\end{align}
where $P_i(\cdot)$ denotes the conditional probability $P_{\cdot|\mathbf{X}}\left(\cdot|\mathbf{x}\in \text{the $i$th class}\right)$ which is just the distribution $N_p(\boldsymbol{\mu}_i,  \boldsymbol{\Sigma})$. The last equality in $\eqref{162}$ is because of the equalities $P_2(A)=P_3(B)$ and $P_1(A)=P_1(B)$ by the symmetry of $\boldsymbol{\mu}_1$, $\boldsymbol{\mu}_2$ and $\boldsymbol{\mu}_3$ and the identity covariance matrix. By the definitions of $T_{OPT}$ and $T$, $\gamma_1$, $\gamma_2$ and $\gamma_3$ are the lines satisfying the equations, $\mathbf{a}_{13}\trans\boldsymbol{\Sigma}^{-1/2}(\mathbf{x}-\mathbf{b}_{13})=0$, $\mathbf{a}_{21}\trans\boldsymbol{\Sigma}^{-1/2}(\mathbf{x}-\mathbf{b}_{21})=0$, and $\mathbf{a}_{23}\trans\boldsymbol{\Sigma}^{-1/2}(\mathbf{x}-\mathbf{b}_{23})=0$, respectively, and $\widetilde{\gamma_1}$, $\widetilde{\gamma_2}$ and $\widetilde{\gamma_3}$ are the lines satisfying equations, $\widehat{\mathbf{a}}_{13}\trans\boldsymbol{\Sigma}^{-1/2}(\mathbf{x}- \widehat{\mathbf{b}}_{13})=0$, $\widehat{\mathbf{a}}_{21}\trans\boldsymbol{\Sigma}^{-1/2}(\mathbf{x}- \widehat{\mathbf{b}}_{21})=0$, and $\widehat{\mathbf{a}}_{23}\trans\boldsymbol{\Sigma}^{-1/2}(\mathbf{x}- \widehat{\mathbf{b}}_{23})=0$, respectively. The following relationship can be found in the proof of Theorem \ref{theorem_1},
\begin{align}
&P \left(\widehat{\mathbf{a}}_{ji}\trans\mathbf{Z}>\widehat{\mathbf{a}}_{ji}\trans\boldsymbol{\Sigma}^{-1/2}(\widehat{\mathbf{b}}_{ji}-\boldsymbol{\mu}_i), \quad \mathbf{a}_{ji}\trans\mathbf{Z}<\mathbf{a}_{ji}\trans\boldsymbol{\Sigma}^{-1/2}(\mathbf{b}_{ji}-\boldsymbol{\mu}_i)\bigg\vert\mathbf{X} \right)\notag\\
=& P_i\left(\widehat{\mathbf{a}}_{ji}\trans\boldsymbol{\Sigma}^{-1/2}(\mathbf{x}- \widehat{\mathbf{b}}_{ji})>0, \quad \mathbf{a}_{ji}\trans\boldsymbol{\Sigma}^{-1/2}(\mathbf{x}-\mathbf{b}_{ji})<0\right).\label{160}
\end{align}  
Therefore, for $i=3$,
\begin{align}
&\sum_{j\neq i}P \left(\widehat{\mathbf{a}}_{ji}\trans\mathbf{Z}>\widehat{\mathbf{a}}_{ji}\trans\boldsymbol{\Sigma}^{-1/2}(\widehat{\mathbf{b}}_{ji}-\boldsymbol{\mu}_i),  \quad\mathbf{a}_{ji}\trans\mathbf{Z}<\mathbf{a}_{ji}\trans\boldsymbol{\Sigma}^{-1/2}(\mathbf{b}_{ji}-\boldsymbol{\mu}_i)\bigg\vert\mathbf{X} \right),\notag\\
=&P_3\left(\widehat{\mathbf{a}}_{13}\trans\boldsymbol{\Sigma}^{-1/2}(\mathbf{x}- \widehat{\mathbf{b}}_{13})>0, \quad \mathbf{a}_{13}\trans\boldsymbol{\Sigma}^{-1/2}(\mathbf{x}-\mathbf{b}_{13})<0\right)\notag\\
&+P_3\left(\widehat{\mathbf{a}}_{23}\trans\boldsymbol{\Sigma}^{-1/2}(\mathbf{x}- \widehat{\mathbf{b}}_{23})>0, \quad \mathbf{a}_{23}\trans\boldsymbol{\Sigma}^{-1/2}(\mathbf{x}-\mathbf{b}_{23})<0\right)\;,\label{166}
\end{align}
where $\{\widehat{\mathbf{a}}_{13}\trans\boldsymbol{\Sigma}^{-1/2}(\mathbf{x}- \widehat{\mathbf{b}}_{13})>0\}$ denotes the half space above the whole line of which $\widetilde{\gamma_1}$ is the left half part, and $\{\mathbf{a}_{13}\trans\boldsymbol{\Sigma}^{-1/2}(\mathbf{x}-\mathbf{b}_{13})<0\}$  denote the half space  below the whole line of which $ {\gamma_1}$ is the left half part. The intersection of these two half spaces is the strip region between the two whole lines and is denoted by $\widetilde{B}$. The event in the second probability on the right hand side of  $\eqref{166}$ is empty and the sum of the two probabilities is $P_3(\widetilde{B})$. Similarly, we can calculate the other two sums for $i=1,2$ and then it can be shown that the right hand side of $\eqref{25}$ is equal to $\frac{2}{3}P_3(\widetilde{B})$. We will show that $P_3(\widetilde{B})/[P_3(B)-P_1(B)]\to 1$. Then in this situation, the upper bound is asymptotically optimal. Let $d_n\to\infty$ and $d_n\epsilon_n\to \infty$. Without loss of generality, we assume that $\gamma_1$ is on the $x$-axis and $\boldsymbol{\mu}_1$ is on the $y$-axis. Then $P_3$ is the two dimensional normal distribution with mean $(0,-d_n)$ and the covariance matrix $\boldsymbol{\Sigma}$ equal to the identity matrix. By the relationship
 \begin{align*}
\{x<\frac{d_n}{2}, -\epsilon_n<y<0\}\subset B \subset \widetilde{B}=\{-\epsilon_n<y<0\},
\end{align*}
we have
\begin{align*}
P_3\left(x<\frac{d_n}{2}, -\epsilon_n<y<0\right)\le P_3(B)\le P_3(\widetilde{B})\le P(-\epsilon_n<y<0),
\end{align*}
and hence 
\begin{align*}
\frac{P_3(\widetilde{B})}{P_3(B)} \le \frac{P_3(-\epsilon_n<y<0)}{P_3\left(x<\frac{d_n}{2}, -\epsilon_n<y<0\right)}=\frac{1}{P_3\left(x<\frac{d_n}{2}\right)}\to 1.
\end{align*}
Similarly, $P_1$ is the two dimensional normal distribution with mean $(0,d_n)$ and the identity covariance matrix. Hence,
\begin{align*}
\frac{P_1(B)}{P_3(\widetilde{B})}\le \frac{P_1(\widetilde{B})}{P_3(\widetilde{B})} \le \frac{P_1(-\epsilon_n<y<0)}{P_3\left(-\epsilon_n<y<0\right)}=\frac{\int_{d_n}^{d_n+\epsilon_n}\phi(x)dx}{\int_{d_n-\epsilon_n}^{d_n}\phi(x)dx}=\frac{\int_{d_n-\epsilon_n}^{d_n}\phi(x+\epsilon_n)dx}{\int_{d_n-\epsilon_n}^{d_n}\phi(x)dx}\le e^{-\frac{d_n\epsilon_n}{2}}\to 0.
\end{align*}
 So we have $P_3(\widetilde{B})/[P_3(B)-P_1(B)]\to 1$. 
 In this example, for simplicity, we make the boundaries of $T$ and $T_{OPT}$ to be parallel, which is not necessary for the bound to be asymptotically optimal. 

We observe that the asymptotic convergence rate based on the upper bound for $R_T/R_{OPT}-1$ obtained when $K>2$ (even if $K$ is fixed) is slower than that when $K=2$. This phenomena can be explained by comparing $\eqref{25}$ and $\eqref{100}$. In the case of $K=2$, there is a negative term on the right hand side of $\eqref{100}$, which has the same order as the positive term and the difference between the two terms has a higher order convergence rate for the classification rules we consider in this paper. But there is no negative terms in $\eqref{25}$. Hence, the asymptotic convergence rate for $K>2$ (even if $K$ is fixed) is slower than that for $K=2$. 
\begin{figure}[h]
\includegraphics[height=3in,width=6in]{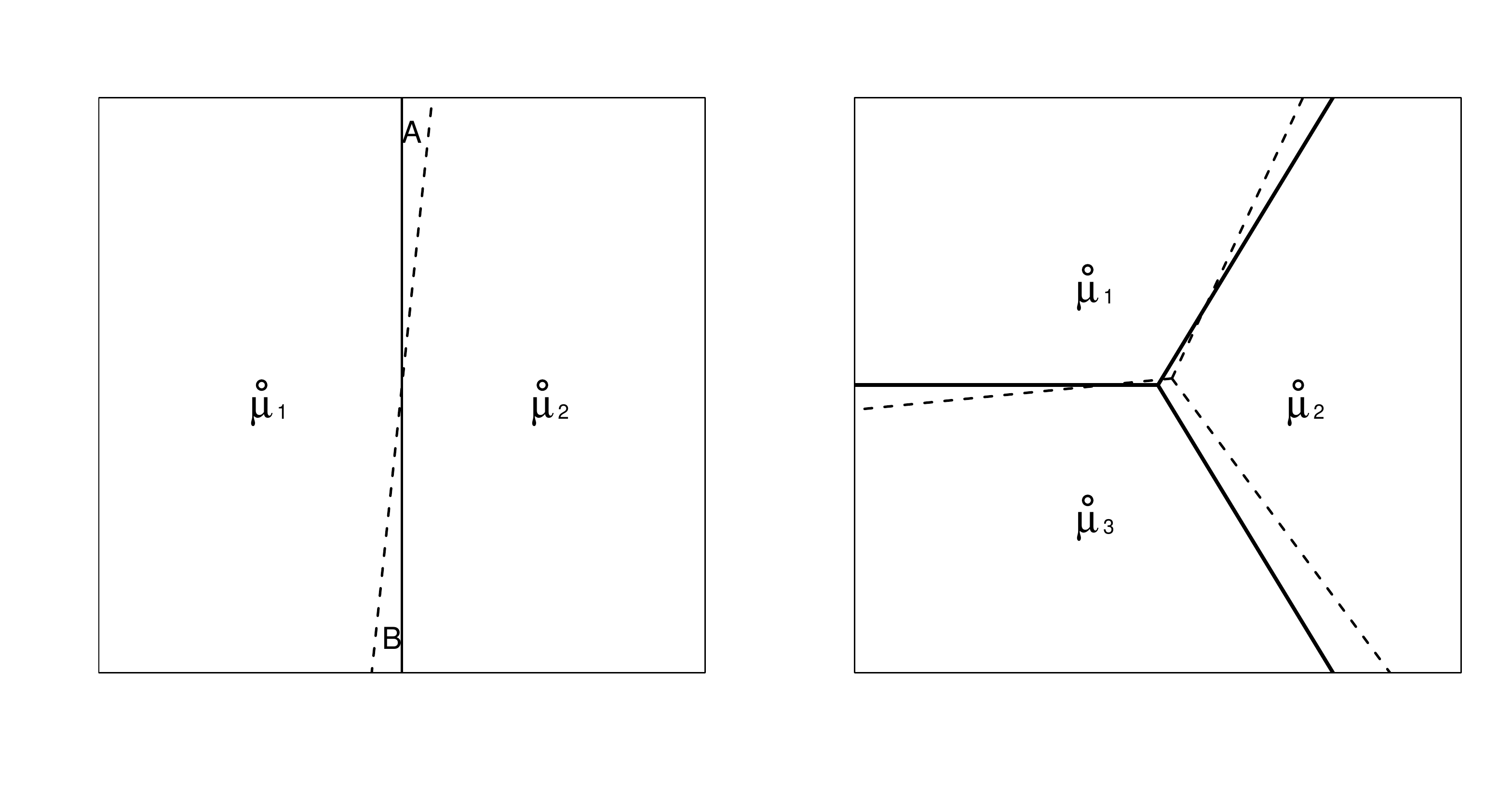}
    \caption{ \baselineskip=10pt The solid lines are the classification boundaries of $T_{OPT}$ and the dashed lines are the classification boundaries of $T$. In the left figure, $K=2$ and $\boldsymbol{\mu}_1$ and $\boldsymbol{\mu}_2$ are the class means. In the right figure, $K=3$.}
\label{fig_1}\end{figure}
We use Figure \ref{fig_1} to illustrate the difference between these two cases. The solid lines are the classification boundaries of $T_{OPT}$ and the dashed lines are the classification boundaries of $T$. In the left figure, $K=2$ and the plot is the projection of the  $p\ge 2$ dimensional feature space  onto the two dimensional space spanned by the vectors orthogonal to the two boundary hyperplanes, respectively. A and B denote the upper and lower regions between the two boundaries. We use $P_1$ and $P_2$ to denote the probability distributions of $N(\boldsymbol{\mu}_1,  \mathbf{I})$ and $N(\boldsymbol{\mu}_2,  \mathbf{I})$, respectively. Using the same arguments as in the above example, we can obtain
\begin{align}
&R_{T}(\mathbf{X})-R_{{OPT}}=\frac{1}{2}[P_1(B)-P_1(A)]+\frac{1}{2}[P_2(A)-P_2(B)]\;.\label{102}
\end{align} 
For the classification rules considered in this paper, A and B actually form a ``matched'' pair in the sense that $P_1(B)-P_1(A)$ is much smaller than $P_1(B)$ and $P_1(A)$, and so is $P_2(A)-P_2(B)$. However, when $K>2$, the boundaries are usually complicated. The right figure is the case of $K=3$ where we cannot find matched pairs as in the case of $K=2$.
 
  Now we consider the asymptotic optimality of a classification rule $T$ in the family.  $\widehat{\mathbf{a}}_{ji}$ and $\widehat{\mathbf{b}}_{ji}$ in the definition $\eqref{24}$ of $T$ are typically estimates of $\mathbf{a}_{ji}$ and $\mathbf{b}_{ji}$ in $\eqref{22}$ for $T_{OPT}$. Given a classification rule $T$ with specific forms of $\widehat{\mathbf{a}}_{ji}$ and $\widehat{\mathbf{b}}_{ji}$, it is relatively easy to calculate the convergence rates of $\widehat{\mathbf{a}}_{ji}$ and $\widehat{\mathbf{b}}_{ji}$ as shown in the following sections.  We will establish the asymptotic optimality of $T$ and the convergence rate for $R_T/R_{OPT}-1$ based on the convergence rates of $\widehat{\mathbf{a}}_{ji}$ and $\widehat{\mathbf{b}}_{ji}$.  		Let $M_{min}=\min_{1\le i\neq j\le K}(\boldsymbol{\mu}_j-\boldsymbol{\mu}_i)\trans\boldsymbol{\Sigma}^{-1}(\boldsymbol{\mu}_j-\boldsymbol{\mu}_i)$ and $M_{max}=\max_{1\le i\neq j\le K}(\boldsymbol{\mu}_j-\boldsymbol{\mu}_i)\trans\boldsymbol{\Sigma}^{-1}(\boldsymbol{\mu}_j-\boldsymbol{\mu}_i)$ be the minimum and maximum  Mahalanobis distances between  any two of the $K$ class means, respectively. By   $\eqref{22}$, $\|\mathbf{a}_{ji}\|^2_2=(\boldsymbol{\mu}_j-\boldsymbol{\mu}_i)\trans\boldsymbol{\Sigma}^{-1}(\boldsymbol{\mu}_j-\boldsymbol{\mu}_i)$. Therefore, under Condition \ref{condition_2}, we have
\begin{align}
c_1\le M_{min}=\min_{i\neq j}\|\mathbf{a}_{ji}\|^2_2\le \max_{i\neq j}\|\mathbf{a}_{ji}\|^2_2=M_{max}.\label{49}
\end{align}
 $M_{min}$ and $M_{max}$ depend on $p$ and $K$, and can go to infinity as $p\to\infty$.

\begin{Theorem}\label{theorem_7}
Suppose that Conditions \ref{condition_1} and \ref{condition_2} hold and $\{s_n: n\ge 1\}$ is a sequence of nonrandom positive numbers with  $M_{max}s_n\to 0$ as $n\to \infty$. For any $1\le j\neq i\le K$, let 
\begin{align*}
\mathbf{a}_{ji}=t_{ji}\widehat{\mathbf{a}}_{ji}+(\mathbf{a}_{ji})_\perp
\end{align*}
be an orthogonal decomposition of $\widehat{\mathbf{a}}_{ji}$, where $t_{ji}\widehat{\mathbf{a}}_{ji}$ is the orthogonal projection of $\mathbf{a}_{ji}$ along the direction of $\widehat{\mathbf{a}}_{ji}$, $t_{ji}$ is a real number, and $(\mathbf{a}_{ji})_\perp$ is orthogonal to $t_{ji}\widehat{\mathbf{a}}_{ji}$.  Let
\begin{align}
&  \widehat{d}_{ji}=\widehat{\mathbf{a}}_{ji}\trans\boldsymbol{\Sigma}^{-1/2}(\widehat{\mathbf{b}}_{ji}-\boldsymbol{\mu}_i)\;, \quad d_{ji}=\mathbf{a}_{ji}\trans\boldsymbol{\Sigma}^{-1/2}(\mathbf{b}_{ji}-\boldsymbol{\mu}_i)=\frac{1}{2}\|\mathbf{a}_{ji}\|_2^2.\notag
\end{align}  
If the following conditions are satisfied,
 \begin{align}
&\|\mathbf{a}_{ji}\|^2_2-\|\widehat{\mathbf{a}}_{ji}\|_2^2=\|\mathbf{a}_{ji}\|_2^2O_p(s_n),\quad t_{ji}=1+O_p(s_n), \quad   d_{ji}-\widehat{d}_{ji}
=\|\widehat{\mathbf{a}}_{ji}\|_2^2O_p(s_n),\label{104}
\end{align} 
where $O_p(s_n)$ are uniform for all $1\le j\neq i\le K$, then we have 
\begin{align}
&\frac{R_{T}(\mathbf{X})}{R_{{OPT}}}-1= O_p\left(K^2\sqrt{M_{max}s_n\log{\left[(M_{max}s_n)^{-1}\right]}} \right).\label{105}
\end{align}
\end{Theorem}
It is natural to ask whether the convergence rate in $\eqref{105}$ can be improved. The following theorem answers the question for the case where $K$ is  bounded.  

\begin{Theorem}\label{theorem_8}
Assume that all the conditions in Theorem \ref{theorem_7} hold. Moreover, suppose that there exists constants $c_4$  and $c_5$ independent of $n$, $p$ and $K$, such that
\begin{align}
&M_{max}/M_{min}\le c_4 ,\quad \min_{1\le i\neq j\le K}\frac{\|(\mathbf{a}_{ji})_\perp\|^2_2}{\|\mathbf{a}_{ji}\|^2_2}>c_5 s_n.\label{186}
\end{align}
 Then with probability converging to 1, the upper bound in $\eqref{25}$  
 \begin{align}
&\frac{1}{K}\sum_{i=1}^K\sum_{j\neq i} P \left(\widehat{\mathbf{a}}_{ji}\trans\mathbf{Z}>\widehat{\mathbf{a}}_{ji}\trans\boldsymbol{\Sigma}^{-1/2}(\widehat{\mathbf{b}}_{ji}-\boldsymbol{\mu}_i),  \mathbf{a}_{ji}\trans\mathbf{Z}<\mathbf{a}_{ji}\trans\boldsymbol{\Sigma}^{-1/2}(\mathbf{b}_{ji}-\boldsymbol{\mu}_i)\bigg\vert\mathbf{X} \right)\notag\\
\ge & \frac{c_5}{4\sqrt{c_4}} \sqrt{M_{max} s_n}R_{OPT}.\label{185}
\end{align}
\end{Theorem}
\begin{Remark}\label{remark_5}
\begin{itemize}
\item[1.] The inequality $\eqref{185}$ indicates that any convergence rate derived from the upper bound in $\eqref{25}$ of Theorem 3.1 is at most  $ \sqrt{M_{max}s_n}$ and cannot be faster. If $K$ is fixed or bounded, the convergence rate in $\eqref{105}$ is different from $ \sqrt{M_{max}s_n}$ only by a multiplying factor $\sqrt{\log{\left[(M_{max}s_n)^{-1}\right]}}$. Hence, the convergence rate in $\eqref{105}$ cannot be improved in this case except for the  logarithm factor as long as the convergence rate is derived from the upper bound in Theorem 3.1.
\item[2.]  The first inequality in $\eqref{186}$ indicates that all the Mahalanobis distances between the $K$ class means have the same orders. If this assumption is not true, then there are some class means  among which the distances will be much smaller than those from them to other class means. $R_{OPT}$ will be dominated by the errors between these classes (here, the error between the $i$th and $j$th classes means the probability that an observation from the $i$th class is assigned to the $j$th class or vice versus).   
\item[3.] The second inequality in $\eqref{186}$ guarantees that the convergence rates of ${\|(\mathbf{a}_{ji})_\perp\|^2_2}/{\|\mathbf{a}_{ji}\|^2_2}$ is slower than $s_n$. The conditions in $\eqref{104}$ in Theorem \ref{theorem_7} implies that the convergence rates of ${\|(\mathbf{a}_{ji})_\perp\|^2_2}/{\|\mathbf{a}_{ji}\|^2_2}$ is $s_n$ or faster. Under these two conditions, the convergence rate of ${\|(\mathbf{a}_{ji})_\perp\|^2_2}/{\|\mathbf{a}_{ji}\|^2_2}$ is exactly $s_n$ in Theorem \ref{theorem_8}. 
\end{itemize}
\end{Remark}

In the following sections, we will apply Theorem \ref{theorem_7} to $T_{LDA}$, the extended $T_{SLDA}$ and $T_{LPD}$.

\section{Classic LDA}
The classification rule $T_{LDA}$ of the classic LDA is a special case of the rule in $\eqref{23}$ with 
\begin{align}
\widehat{\mathbf{a}}_{ji}= {\boldsymbol{\Sigma}}^{1/2}\widehat{\boldsymbol{\Sigma}}^{-1}(\bar{\mathbf{x}}_j-\bar{\mathbf{x}}_i), \quad \widehat{\mathbf{b}}_{ji}=(\bar{\mathbf{x}}_i+\bar{\mathbf{x}}_j)/2,\label{110}
\end{align}
 for all $1\le i\neq j\le K$.  Let $\widehat{\sigma}_{kl}$ and $\sigma_{kl}$ denote the $(k,l)$-th elements of $\widehat{\boldsymbol{\Sigma}}$ and $\boldsymbol{\Sigma}$, respectively, where $1\le k,l\le p$. When there is one population, that is, $K=1$, it has been shown in (12) of \citet{bickel2008covariance} that $\max_{k,l} \left|\widehat{\sigma}_{kl}-\sigma_{kl}\right|=O_p\left(\sqrt{ \log{p}/{n}}\right)$. Because we allow the number of classes $K$ to go to infinity, in the following theorem, we will consider the effect of the large $K$.
\begin{Theorem}\label{theorem_3}
Suppose that Condition \ref{condition_1} holds. Then we have
\begin{align*}
&\max_{k,l} \left|\widehat{\sigma}_{kl}-(1-\frac{K}{n})\sigma_{kl}\right|=O_p\left(\sqrt{\frac{\log{p}}{n}}\right). 
\end{align*}
Moreover, 
\begin{align}
&\|\widehat{\boldsymbol{\Sigma}}-(1-\frac{K}{n})\boldsymbol{\Sigma}\|=O_p\left(p\sqrt{\frac{\log{p}}{n}}\right), \text{ and hence } \|\widehat{\boldsymbol{\Sigma}}-\boldsymbol{\Sigma}\|=O_p\left(p\sqrt{\frac{\log{p}}{n}}+\frac{K}{n}\right).\label{34}
\end{align}
\end{Theorem}
 From Theorem \ref{theorem_3}, one can see that a large $K$ has a shrinkage effect on $\widehat{\boldsymbol{\Sigma}}$, that is, when $K$ is large, the entry $\widehat{\sigma}_{kl}$ of $\widehat{\boldsymbol{\Sigma}}$ is close to the shrunk entry of $\boldsymbol{\Sigma}$.
The convergence rates in Theorem \ref{theorem_3} plays a basic role  in the following theoretical development.

\citet{Shao-2011} provided necessary conditions for the classic LDA to be asymptotically optimal in the case $K=2$ as both $n, p\to \infty$. We will extend the results to the case of $K>2$. Before we state the theorem, we will exclude the situations where there are very small numbers of observations in some classes by assuming the following condition:
\begin{Condition}\label{condition_3}
There exists a constant $c_2$ independent of $n$, $p$, and $K$ such that 
\begin{align}
\frac{1}{\min_{1\le i\le K} n_i}\le c_2\left(\frac{K}{n}\right).\label{120}
\end{align}
\end{Condition}
This condition implies that $n_j\ge \min_{1\le i\le K} n_i\ge c_2^{-1}n/K$ for any $1\le j\le K$. Therefore, the number of observations in any class is of the same order as the average number $n/K$ of observations. 
\begin{Theorem}\label{theorem_2}
 Suppose that Conditions \ref{condition_1}-\ref{condition_3} hold and $K\le p+1$. Let 
\begin{align}
s_n=p\sqrt{\frac{\log{p}}{n}}+\frac{K}{\sqrt{M_{min}}}\sqrt{\frac{p}{n}}\to 0. \label{125}
\end{align}
 Then if $K^2\sqrt{M_{max}s_n\log{\left[(M_{max}s_n)^{-1}\right]}}\to 0$, the classic LDA is asymptotically optimal and 
\begin{align}
&\frac{R_{LDA}(\mathbf{X})}{R_{{OPT}}}-1 = O_p\left(K^2\sqrt{M_{max}s_n\log{\left[(M_{max}s_n)^{-1}\right]}} \right). \label{128}
\end{align}
\end{Theorem}
\begin{Remark}\label{remark_3}
\begin{itemize}
\item[1.] We compare Theorem \ref{theorem_2} with Theorem 1 in \citet{Shao-2011} where it is assumed that only the first term  $p\sqrt{\log{p}/n}$ of $s_n$ in $\eqref{125}$ converges to zero. The second term of $s_n$ in $\eqref{125}$ is due to the effect of $K$. When $K$ has the order $\sqrt{M_{min}p\log{p}}$ or larger, the second term in $s_n$ is not negligible.  
\item[2.] If $K$ is fixed or bounded and $p\to \infty$, then the second term of $s_n$ satisfies
\begin{align*}
\frac{K}{\sqrt{M_{min}}}\sqrt{\frac{p}{n}}=o\left(p\sqrt{\frac{\log{p}}{n}}\right).
\end{align*}
In this case, the condition $K^2\sqrt{M_{max}s_n\log{\left[(M_{max}s_n)^{-1}\right]}}\to 0$ in Theorem \ref{theorem_2} is equivalent to $M_{max}s_n \to 0$. When $K=2$, this is the same as that  in Theorem 1 in \citet{Shao-2011}. To see this, note that when $K=2$,
\[M_{min}=(\boldsymbol{\mu}_1-\boldsymbol{\mu}_2)\trans\boldsymbol{\Sigma}^{-1}(\boldsymbol{\mu}_1-\boldsymbol{\mu}_2)= M_{max}\;,\]
 and hence the $\Delta_p^2$ defined as $(\boldsymbol{\mu}_1-\boldsymbol{\mu}_2)\trans\boldsymbol{\Sigma}^{-1}(\boldsymbol{\mu}_1-\boldsymbol{\mu}_2)$ in \citet{Shao-2011} is equal to  $ M_{max}$. Therefore, $M_{max}s_n \to 0$ is equivalent to $\Delta_p^2s_n\to 0$ which is the condition in Theorem 1 in \citet{Shao-2011}. Moreover, as discussed in Section 3, when $K>2$, the convergence rate in $\eqref{128}$ is smaller than that for $K=2$ in Theorem 1 in \citet{Shao-2011}. 
\end{itemize}
\end{Remark}

\section{Sparse LDA by thresholding}\label{sect4}
It has been shown that when $p/n\to \infty$, the classic LDA may not be asymptotically optimal in Theorem 2 of \citet{Shao-2011}. By imposing sparsity conditions on $\boldsymbol{\Sigma}$ and $\boldsymbol{\mu}_1-\boldsymbol{\mu}_2$, \citet{Shao-2011} proposed a sparse LDA rule by thresholding and proved that it is asymptotically optimal in the case of $K=2$. In this section, we extend this method to arbitrary $K$ and provide asymptotic optimality.. As in \citet{Shao-2011}, we consider the following sparsity measure on $\boldsymbol{\Sigma}$ in \citet{bickel2008covariance},
\begin{align}
C_{h,p}=\max_{1\le k\le p}\sum_{l=1}^p|\sigma_{kl}|^h,\label{51}
\end{align}
where $0\le h<1$ is a constant independent of $p$. \citet{Shao-2011} used the sparse estimate of $\boldsymbol{\Sigma}$ in \citet{bickel2008covariance} by performing a thresholding procedure on $\widehat{\boldsymbol{\Sigma}}$. In tis case, we need to consider the effect of large $K$.  Hence, we propose to apply the thresholding procedure  to $(1-K/n)^{-1}\widehat{\boldsymbol{\Sigma}}$, instead of to $\widehat{\boldsymbol{\Sigma}}$. Specifically, let $\widetilde{\boldsymbol{\Sigma}}$ be the thresholding estimate with the $(k,l)$th entry 
\begin{align*}
\widetilde{\sigma}_{kl}=(1-K/n)^{-1}\widehat{\sigma}_{kl}I_{[(1-K/n)^{-1}|\widehat{\sigma}_{kl}|\ge t_n]}, \forall 1\le k,l\le p,
\end{align*}
where $t_n=M_1\sqrt{\frac{\log{p}}{n}}$, $M_1$ is a large enough positive constant and $\widehat{\sigma}_{kl}$ is the  $(k,l)$th entry of $\widehat{\Sigma}$.  We first derive the convergence rate for the revised thresholding estimator.

\begin{Theorem}\label{theorem_4}
Suppose that Condition \ref{condition_1} holds, $\log{p}/n=o(1)$ and
 \begin{align*}
d_n=C_{h,p}\left(\frac{\log{p}}{n}\right)^{(1-h)/2}=o(1).
\end{align*}
Then
 \begin{align}
\|\widetilde{\boldsymbol{\Sigma}}-\boldsymbol{\Sigma}\|=O_p(d_n),\quad \|\widetilde{\boldsymbol{\Sigma}}^{-1}-\boldsymbol{\Sigma}^{-1}\|=O_p(d_n). \label{64}
\end{align}
\end{Theorem}
  $d_n$ is  the convergence rate of the thresholding estimate of the covariance matrix for one population in \citet{bickel2008covariance}. Hence, the revised thresholding estimate has the same convergence rates for any $K$. Now define 
\begin{align}
\boldsymbol{\delta}_{ji}=\boldsymbol{\mu}_j-\boldsymbol{\mu}_i,\quad \widehat{\boldsymbol{\delta}}_{ji}=\bar{\mathbf{x}}_j-\bar{\mathbf{x}}_i\label{130}
\end{align}
 for all $1\le i, j\le K$. We define the following sparsity measure on $\boldsymbol{\delta}_{ji}$ which is an extension of (9) in \citet{Shao-2011},
 \begin{align}
D_{g,p}=\max_{1\le i\neq j\le K}\sum_{k=1}^p(\boldsymbol{\delta}_{ji}^k)^{2g}, \label{68}
\end{align}
where $\boldsymbol{\delta}_{ji}^k$ is the $k$th coordinate of $\boldsymbol{\delta}_{ji}$ and $0\le g<1$ is a constant independent of $p$ and $K$. When $K=2$, there is essentially one $\boldsymbol{\delta}_{ji}$ because $\boldsymbol{\delta}_{21}=-\boldsymbol{\delta}_{12}$. In this case, the $D_{g,p}$ is just that in \citet{Shao-2011}. We extend the classification rule of the sparse LDA method in \citet{Shao-2011} for arbitrary $K$ as follows.
\begin{align*}
T_{SLDA} &\text{: to allocate a new observation $\mathbf{x}$ to the $i$ class if }  \widetilde{\boldsymbol{\delta}}_{ji}\trans\widetilde{\boldsymbol{\Sigma}}^{-1}(\mathbf{x}- \widetilde{\mathbf{b}}_{ji})<0,\quad \forall j\neq i,
\end{align*}
where $\widetilde{\boldsymbol{\delta}}_{ji}$ is a sparse estimator of $\boldsymbol{\delta}_{ji}$ and $\widetilde{\mathbf{b}}_{ji}$ is an estimate of $(\boldsymbol{\mu}_j+\boldsymbol{\mu}_i)/2$. One may naturally take $\widetilde{\boldsymbol{\delta}}_{ji}$ to be $\widehat{\boldsymbol{\delta}}_{ji}$ thresholded and $\widetilde{\mathbf{b}}_{ji}=(\bar{\mathbf{x}}_i+\bar{\mathbf{x}}_j)/2$ as in \citet{Shao-2011}. However, this choice of $\widetilde{\boldsymbol{\delta}}_{ji}$ and $\widetilde{\mathbf{b}}_{ji}$ is problematic because in this case, there may exist multiple $i$'s which satisfy $\widetilde{\boldsymbol{\delta}}_{ji}\trans\widetilde{\boldsymbol{\Sigma}}^{-1}(\mathbf{x}- \widetilde{\mathbf{b}}_{ji})<0$, for all $j\neq i$, and hence we cannot uniquely determine the class to which we assign $\mathbf{x}$. Therefore, we propose the following estimates. Define the following threshold for $\boldsymbol{\delta}_{ji}$,
\begin{align}
a_n=M_2\left(\frac{\log{p}}{n}\right)^{\alpha} \label{69}
\end{align}
with $M_2>0$  and $\alpha\in (0,1/2)$ are constants.  Define the thresholding estimates
\begin{align*}
 \widetilde{\boldsymbol{\delta}}_{11}=\mathbf{0}, \quad \widetilde{\boldsymbol{\delta}}_{j1}=\widehat{\boldsymbol{\delta}}_{j1}I_{(|\widehat{\boldsymbol{\delta}}_{j1}|\ge a_n)}, \quad  1<j\le K,   
\end{align*}
and let
\begin{align*}
  \widetilde{\boldsymbol{\delta}}_{ji}=\widetilde{\boldsymbol{\delta}}_{j1}-\widetilde{\boldsymbol{\delta}}_{i1}, \quad \widetilde{\mathbf{b}}_{ji}=\bar{\mathbf{x}}_1+\frac{\widetilde{\boldsymbol{\delta}}_{j1}+\widetilde{\boldsymbol{\delta}}_{i1}}{2}, \quad \forall 1\le i, j\le K,
\end{align*}
where for simplicity, we first estimate $\widetilde{\boldsymbol{\delta}}_{j1}$, $1\le j\le K$, but one can first estimate $\widetilde{\boldsymbol{\delta}}_{j2}$, $1\le j\le K$, without any effects on the asymptotic optimality results. Under the above definitions, $\widetilde{\boldsymbol{\delta}}_{ji}\trans\widetilde{\boldsymbol{\Sigma}}^{-1}(\mathbf{x}- \widetilde{\mathbf{b}}_{ji})<0$ if and only if
\begin{align*}
\left(\mathbf{x}-(\widetilde{\boldsymbol{\delta}}_{i1}+\bar{\mathbf{x}}_1)\right)\trans\widetilde{\boldsymbol{\Sigma}}^{-1}\left(\mathbf{x}-(\widetilde{\boldsymbol{\delta}}_{i1}+\bar{\mathbf{x}}_1)\right)<\left(\mathbf{x}-(\widetilde{\boldsymbol{\delta}}_{j1}+\bar{\mathbf{x}}_1)\right)\trans\widetilde{\boldsymbol{\Sigma}}^{-1}\left(\mathbf{x}-(\widetilde{\boldsymbol{\delta}}_{j1}+\bar{\mathbf{x}}_1)\right).
\end{align*}
Hence, $T_{SLDA}$ assigns $\mathbf{x}$ to the $i$th class if and only if $\left(\mathbf{x}-(\widetilde{\boldsymbol{\delta}}_{i1}+\bar{\mathbf{x}}_1)\right)\trans\widetilde{\boldsymbol{\Sigma}}^{-1}\left(\mathbf{x}-(\widetilde{\boldsymbol{\delta}}_{i1}+\bar{\mathbf{x}}_1)\right)$ is the smallest among $\left(\mathbf{x}-(\widetilde{\boldsymbol{\delta}}_{j1}+\bar{\mathbf{x}}_1)\right)\trans\widetilde{\boldsymbol{\Sigma}}^{-1}\left(\mathbf{x}-(\widetilde{\boldsymbol{\delta}}_{j1}+\bar{\mathbf{x}}_1)\right)$, $1\le j\le K$. It is easy to see that $T_{SLDA}$ is a special case of $\eqref{23}$ with
\begin{align}
\widehat{\mathbf{a}}_{ji}= {\boldsymbol{\Sigma}}^{1/2}\widetilde{\boldsymbol{\Sigma}}^{-1}\widetilde{\boldsymbol{\delta}}_{ji}, \quad \widehat{\mathbf{b}}_{ji}=\widetilde{\mathbf{b}}_{ji}=\bar{\mathbf{x}}_1+\frac{\widetilde{\boldsymbol{\delta}}_{j1}+\widetilde{\boldsymbol{\delta}}_{i1}}{2},\quad \forall 1\le i\neq j\le K.\notag
\end{align} 
Moreover, Let $r>1$ be a fixed constant and define 
\begin{align*}
q_n=\max_{1<j\le K}\{\text{the number of $k$'s with  $|\boldsymbol{\delta}_{j1}^k|>a_n/r$}\}.
\end{align*}
Then by Lemma 2 (ii) in \citet{Shao-2011}, with probability converging to 1, the number of nonzero coordinates of $\widetilde{\boldsymbol{\delta}}_{j1}$ is less than or equal to $q_n$ and hence, the number of nonzero coordinates of $\widetilde{\boldsymbol{\delta}}_{ji}$ is less than or equal to $2q_n$.

\begin{Theorem}\label{theorem_5}
 Suppose that Conditions \ref{condition_1}-\ref{condition_3} hold and
\begin{align}
b_n= \max\left\{d_n, \frac{\sqrt{a_n^{2(1-g)}D_{g,p} }}{
\sqrt{M_{min}}},\frac{\sqrt{{(C_{h,p}+K)q_n}}}
{\sqrt{nM_{min}}}\right\}\to 0,\label{82}
\end{align}
where $d_n$ is defined in Theorem \ref{theorem_4}. Then $T_{SLDA}$ is asymptotically optimal and 
\begin{align}
&\frac{R_{SLDA}(\mathbf{X})}{R_{{OPT}}}-1 = O_p\left(K^2\sqrt{M_{max}b_n\log{\left[(M_{max}b_n)^{-1}\right]}} \right) \label{10001}
\end{align}
\end{Theorem}
 The only difference between $b_n$ in $\eqref{82}$ and that in Theorem 3 in \citet{Shao-2011} is that we use $C_{h,p}+K$ in $\eqref{82}$ to replace $C_{h,p}$ in \citet{Shao-2011}.  

\section{Linear programming discriminant rule}
\citet{cai2011direct} observed that  the optimal classification rule $T_{OPT}$ depends on ${\boldsymbol{\Sigma}}$ only through the vectors $\boldsymbol{\beta}_{ji}={\boldsymbol{\Sigma}}^{-1}\boldsymbol{\delta}_{ji}$,  where $\boldsymbol{\delta}_{ji}$ is defined in $\eqref{130}$, $1\le i,j\le K$. They proposed the linear programming discriminant rule where the key step is to estimate $\boldsymbol{\beta}$ through a constrained $l_1$ minimization. We extend the LPD rule to the case where $K$ is arbitrary.  We first define the estimates $\widehat{\boldsymbol{\beta}}_{11}=\mathbf{0}$ and $\widehat{\boldsymbol{\beta}}_{j1}$ to be the solution to 
\begin{align}
\min_{\boldsymbol{\beta}\in \mathbb{R}^p}\|\boldsymbol{\beta}\|_1,\quad \text{subject to }\|\overline{\boldsymbol{\Sigma}}\boldsymbol{\beta}-\widehat{\boldsymbol{\delta}}_{j1}\|_\infty\le \lambda_n,\label{137}
\end{align}
for any $1<j\le K$, where $\overline{\boldsymbol{\Sigma}}=(1-K/n)^{-1}\widehat{\boldsymbol{\Sigma}},\quad \lambda_n=C\sqrt{M_{max}\log{p}/n}$. Then let
\begin{align*}
\widehat{\boldsymbol{\beta}}_{ji}=\widehat{\boldsymbol{\beta}}_{j1}-\widehat{\boldsymbol{\beta}}_{i1},
\end{align*}
for any $1\le i,j \le K$. We use $\overline{\boldsymbol{\Sigma}}$ in the optimization problem $\eqref{137}$ instead of $\widehat{\boldsymbol{\Sigma}}$ as in \citet{cai2011direct} to remove the shrinkage effect of large $K$. By Theorem \ref{theorem_3},
\begin{align}
\|\overline{\boldsymbol{\Sigma}}-\boldsymbol{\Sigma}\|_{\infty}=O_p\left(\sqrt{\frac{\log{p}}{n}}\right).\label{142}
\end{align}
  The classification rule is
\begin{align*}
T_{LPD} &\text{: to allocate a new observation $\mathbf{x}$ to the $i$th class if }  \widehat{\boldsymbol{\beta}}_{ji}\trans (\mathbf{x}- \frac{\bar{\mathbf{x}}_j+\bar{\mathbf{x}}_i}{2})<0, \quad \forall j\neq i,
\end{align*}
which is of the general form $\eqref{23}$ with
\begin{align*}
\widehat{\mathbf{a}}_{ji}={\boldsymbol{\Sigma}}^{1/2}\widehat{\boldsymbol{\beta}}_{ji},\quad \widehat{\mathbf{b}}_{ji}=\frac{\bar{\mathbf{x}}_j+\bar{\mathbf{x}}_i}{2}. 
\end{align*}

\begin{Theorem}\label{theorem_6}
 Suppose that Conditions \ref{condition_1}-\ref{condition_3} hold and
 \begin{align}
r_n=  \left(\sqrt{KM_{max}}\max_{1\le i\neq j\le K}\|\boldsymbol{\Sigma}^{-1}\boldsymbol{\delta}_{ji}\|_1+\max_{1\le i\neq j\le K}\|\boldsymbol{\Sigma}^{-1}\boldsymbol{\delta}_{ji}\|_1^2\right)\sqrt{\frac{\log{p}}{n}}\to 0.\label{90}
\end{align}
Then $T_{LPD}$ is asymptotically optimal and 
\begin{align}
&\frac{R_{LPD}(\mathbf{X})}{R_{{OPT}}}-1 = O_p\left(K^2\sqrt{M_{max}r_n\log{\left[(M_{max}r_n)^{-1}\right]}} \right). \label{10002}
\end{align}
\end{Theorem}
  If $K$ is fixed or bounded. The condition is essentially the same as that in Theorem 3 of \citet{cai2011direct}.

\section{Simulation Study}
In this section, we conduct simulation studies to evaluate the  classification performance of the extended LPD and SLDA,  and compare them with  the nearest shrunken centroids method (NSC), the classic LDA rule with a generalized inverse matrix (GLDA) and the optimal classification rule. 
There are several different definitions for the generalized inverse matrix. We use the Moore-Penrose pseudoinverse and calculate it using the matlab function ``{\it pinv}''.  Although we cannot calculate the optimal rule in practice, we include it in this simulation study as a benchmark. All the methods are implemented in matlab except NSC, which is implemented using the R package ``{\it pamr}'' with the default setting and cross-validation procedure.  

 In LPD, we choose the tuning parameter $\lambda_n$ from the set $\{0.2, 0.25, 0.3, \ldots, 0.65, 0.7\}$ in $\eqref{137}$ by the five-fold cross-validation procedure. There are two tuning parameters in SLDA: $t_n=M_1\sqrt{\log{p}/n}$ and $a_n=M_2\left(\log{p}/n\right)^{\alpha}$, the thresholds for  $\boldsymbol{\Sigma}$ and $\boldsymbol{\delta}_{ji}$, respectively. $M_1$ is chosen from  $\{10^{-5}, 10^{-4}, \ldots, 1\}$  and $M_2$ is from $\{10^{-7}, 10^{-6}, \ldots, 1\}$. We choose $\alpha=0.3$ as in \citet{Shao-2011}. One practical issue of SLDA is that the thresholded covariance matrix may not  be invertible. We propose two different ways to overcome this problem and compare their performance in this study. In the first one,  we use the Moore-Penrose pseudoinverse  of $\widetilde{\boldsymbol{\Sigma}}$ to replace $\widetilde{\boldsymbol{\Sigma}}^{-1}$ in $T_{SLDA}$ and denote this approach by SLDA1. In the second one, we   replace $\widetilde{\boldsymbol{\Sigma}}^{-1}$ by $(\widetilde{\boldsymbol{\Sigma}}+\varepsilon\mathbf{I}_p)^{-1}$, where $\varepsilon$ is a small positive number which will be viewed as a tuning parameter and $\mathbf{I}_p$ is the $p$-dimensional identity matrix. $(\widetilde{\boldsymbol{\Sigma}}+\varepsilon\mathbf{I}_p)^{-1}$ is the inverse of $\widetilde{\boldsymbol{\Sigma}}+\varepsilon\mathbf{I}_p$ if it is full rank, otherwise, $(\widetilde{\boldsymbol{\Sigma}}+\varepsilon\mathbf{I}_p)^{-1}$ is the generalized inverse. This approach is denoted by SLDA2. We choose $\varepsilon$ from the set $\{10^{-5}, 10^{-4}, 10^{-3},10^{-2},10^{-1}\}$. There are two tuning parameters for SLDA1 and three for SLDA2. For both SLDA1 and SLDA2, we choose the tuning parameters by the five-fold cross-validation procedure.  

 We will consider three models with different class means and within-class covariance matrices. For each model, we  consider three different numbers of classes $K=3,6, 9$ and two different dimensionality $p=300, 600$. The total number of observations from all the classes is fixed to be $n=450$ for all models. The numbers of observations from all the classes are the same, that is, $n_1=n_2=\cdots=n_K=450/K$. Therefore, when $K=3$, we have 150 observations from each class, and when $K=9$, we have only 50 observations from each class. For all the three models, the class means have the forms: $\boldsymbol{\mu}_k=(\underbrace{0, \ldots, 0}_{\substack{(k-1)s_0}}, \underbrace{1, \ldots, 1}_{\substack{s_0}}, \underbrace{0, \ldots, 0}_{\substack{p-ks_0}})'$, $1 \le k \le K$. There are $s_0$ ones and all the other numbers are equal to zero in the $p$-dimensional vectors of the class means. The specific details of the three models are given below:
\begin{itemize}
\item[1.] 
$s_0=5 $ and $\boldsymbol{\Sigma}=(\sigma_{ij})_{p \times p}$ with  $\sigma_{ii}=1$ for $1 \le i \le p$ and $\sigma_{ij}=0.5$ for $i \ne j$.
\item[2.] 
$s_0=3$ and $\boldsymbol{\Sigma}$ is a diagonal block matrix  given by
\begin{align}
\boldsymbol{\Sigma} = \begin{bmatrix}
      \boldsymbol{\Sigma}_{11} &  \mathbf{0}           \\ 
      \mathbf{0}          & \boldsymbol{\Sigma}_{22}\\ 
         \end{bmatrix},
\end{align}
where $\boldsymbol{\Sigma}_{11}$ is a $100\times 100$ matrix with  diagonal element equal to 1 and off-diagonal element equal to 0.7, and $\boldsymbol{\Sigma}_{22}$ is a $(p-100)\times (p-100)$ matrix with  diagonal element equal to 1 and off-diagonal element equal to 0.5.
\item[3.]
$s_0=10$ and  $\boldsymbol{\Sigma}=(\sigma_{ij})_{p \times p}$ with  $\sigma_{ij}=0.95^{|i-j|}$ for $i \ne j$.
\end{itemize}
Models 1 and 3 are similar to those in \cite{cai2011direct}. In Model 1,  all the entries of ${\boldsymbol{\Sigma}}^{-1}$ and ${\boldsymbol{\Sigma}}^{-1} \delta_{ij}$ ($i \ne j$) are nonzero. In Model 2,  ${\boldsymbol{\Sigma}}^{-1}$ is also a diagonal block matrix with two diagonal blocks equal to  ${\boldsymbol{\Sigma}}_{11}^{-1}$ and ${\boldsymbol{\Sigma}}_{22}^{-1}$, respectively.  All the entries of ${\boldsymbol{\Sigma}}_{11}^{-1}$ and ${\boldsymbol{\Sigma}}_{22}^{-1}$ are nonzero, but all the coordinates of the vector ${\boldsymbol{\Sigma}}^{-1} \delta_{ij}$ ($i \ne j$)  are equal to zero except the first 100 entries.   In Model 3, the entries of $\boldsymbol{\Sigma}$ decays exponentially with $|i-j|$ and hence satisfies the conditions for SLDA. The $(i,j)$-th entry of $\boldsymbol{\Sigma}^{-1}$ is zero when $|i-j|>1$. Therefore, $\boldsymbol{\Sigma}^{-1}$ is a  sparse matrix. 

For each setting in each model, we repeat the following procedure 50 times. In each repeat, we generate $n=450$ training samples with $n/K$ samples in each class, and then generate 450 test samples with $n/K$ samples in each class independent of the training samples. For each method (except the optimal rule), we use the training data to choose the tuning parameters and find the classification rule, and then we apply the fitted classification rule to the test data to obtain the classification errors.  The average classification errors and standard deviations of the  50 replications for all the setting in the three models are listed in Table \ref{sim2.tab}. 

In all settings, SLDA2 has lower classification errors than SLDA1. Therefore, for SLDA, first adding a diagonal matrix with small common positive diagonal entries to the thresholded  covariance matrix is better than directly calculating its generalized inverse matrix in terms of predictive performance. For all settings in both Models 1 and 2, LPD has the smallest classification errors compared to SLDA1, SLDA2, NSC and GLDA. In Model 3, when $p=300$, LPD has the smallest errors and when $p=600$, SLDA2 has the smallest errors.  For all the three models and all the methods, the classification errors increase with the increasing number of classes when $p$ is fixed. Given $K$, the optimal errors are almost unchanged when $p$ increases from 300 to 600, and the classification errors of all other methods increase.  When both $K$ and $p$ are large,  GLDA performs much worse than other methods and have errors close to or more than 50\%.  

 To evaluate the computational efficiency of the extended methods: LPD, SLDA1 and SLDA2, we list the  average time in seconds of running one replication (including 5 fold cross-validation) for Model 1 of the three methods in Table \ref{time.tab}. All the computations are conducted on a compute cluster with Red Hat Linux. The CPU in each node of the cluster is {\it : Intel(R) Xeon(R) CPU  5160  \@ 3.00GHz}. All methods need longer execution time as  $p$ increases with $K$ fixed. Compared to SLDA1 and SLDA2, the execution time of LPD  is much more sensitive to the number of classes. That is because the major computation load of LPD is to solve the optimization problem $\eqref{137}$ for each $1<j\le K$, and the major load of SLDA is the calculation of the inverse or generalized inverse matrix of the thresholded covariance matrix with or without a small diagonal matrix .

\begin{table}[h]
\caption{\baselineskip=10pt The averages and standard deviations (in parenthesis) of classification errors of 50 replications for all the settings in all the three models.  }
\vspace{5mm}
\label{sim2.tab}\centering
 \small\addtolength{\tabcolsep}{-3pt}
 {
\begin{tabular}{|c|c|c|c|c|c|c|c|}
\hline
$K$ & $p$ & LPD & SLDA1 & SLDA2 & NSC & GLDA & Optimal \\\hline
\multicolumn{8}{|c|}{Model 1}\\\hline
\multirow{2}{*}{3} & 300 &
0.027(0.007)& 0.171(0.027)& 0.067(0.013)& 0.052(0.026)& 0.197(0.023)& 0.023(0.006)\\\cline{2-8}
&{600}&
0.030(0.008)& 0.273(0.032)& 0.083(0.012)& 0.052(0.031)& 0.312(0.034)& 0.023(0.007)\\\hline  
\multirow{2}{*}{6} & {300}&
0.063(0.012)& 0.314(0.038)& 0.150(0.032)& 0.100(0.042)& 0.379(0.032)& 0.050(0.011)\\\cline{2-8}
&{600}&
0.079(0.016)& 0.456(0.041)& 0.167(0.022)& 0.091(0.039)& 0.543(0.029)& 0.048(0.011)\\\hline
\multirow{2}{*}{9} & {300}&
0.095(0.015)& 0.414(0.040)& 0.211(0.037)& 0.140(0.035)&  0.495(0.029)& 0.071(0.013)\\\cline{2-8}
&{600}&
0.134(0.020)& 0.552(0.052)& 0.251(0.047)& 0.152(0.051)& 0.659(0.024)& 0.073(0.013)\\\hline
 \multicolumn{8}{|c|}{Model 2}\\\hline
\multirow{2}{*}{3} & {300} &
0.027(0.008)& 0.176(0.025)& 0.066(0.016)& 0.068(0.042)& 0.198(0.025)& 0.024(0.008) \\\cline{2-8}
&{600}&
0.027(0.007)& 0.313(0.028)& 0.083(0.013)& 0.068(0.048)& 0.373(0.026)& 0.023(0.008)\\\hline
\multirow{2}{*}{6} & {300}&
0.069(0.035)& 0.324(0.040)& 0.144(0.036)& 0.109(0.045)&  0.362(0.059)& 0.047(0.016)\\\cline{2-8}
&{600}&
0.066(0.014)& 0.509(0.043)& 0.189(0.053)& 0.109(0.046)& 0.597(0.026)& 0.048 (0.009)\\\hline
\multirow{2}{*}{9} & {300} &
0.091(0.013)& 0.433(0.047)& 0.231(0.055)& 0.169(0.056)&  0.502(0.022)& 0.072(0.012)\\\cline{2-8}
&{600}&
0.110(0.017)& 0.618(0.055)& 0.292(0.077)& 0.181(0.047)&  0.711(0.024)& 0.073 (0.013)\\\hline
\multicolumn{8}{|c|}{Model 3}\\\hline
\multirow{2}{*}{3} & {300} &
0.006(0.004)& 0.066(0.016)& 0.020(0.007)& 0.190(0.029)&  0.068(0.017)& 0.002(0.002)\\\cline{2-8}
&{600}&
0.084(0.039)& 0.673(0.049)& 0.034(0.010)& 0.187(0.032)& 0.160(0.017)& 0.002(0.002)\\\hline
\multirow{2}{*}{6} & {300} &
0.016(0.006)& 0.171(0.027)& 0.057(0.013)& 0.384(0.025)& 0.176(0.024)& 0.004(0.003)\\\cline{2-8}
&{600}&
0.111(0.062)& 0.835(0.028)& 0.091(0.024)& 0.382(0.025)& 0.337(0.025)& 0.006(0.005)\\\hline
\multirow{2}{*}{9} & {300} &
0.031(0.008)& 0.262(0.029)& 0.098(0.016)& 0.501(0.030)& 0.257(0.024)& 0.008(0.004)\\\cline{2-8}
&{600}&
0.185(0.069)& 0.897(0.019)& 0.170(0.046)& 0.508(0.024)& 0.456(0.028)& 0.007(0.005)\\\hline
 \end{tabular}}
\end{table}

\begin{table}[h]
\caption{\baselineskip=10pt The  average running time (in seconds) and standard deviations (in parenthesis) of one replication (including 5 fold cross-validation) for Model 1. }
\vspace{5mm}
\label{time.tab}\centering
 \small\addtolength{\tabcolsep}{-3pt}
 {
\begin{tabular}{|c|c|c|c|c|}
\hline
$K$ & $p$ & LPD & SLDA1 & SLDA2   \\\hline
\multirow{2}{*}{3} &  
300& 9.29(1.12)& 19.15(0.64)& 96.84(2.83) \\ 
&600& 59.42(3.22)& 58.11(5.84)& 309.68(9.57)\\ 
\hline
 \multirow{2}{*}{6} &  
300& 25.77(2.35)& 22.28(3.45)& 113.98(16.91)\\ 
&600& 157.52(2.44)& 57.53(4.11) & 310.29(7.75) \\ 
\hline
\multirow{2}{*}{9} &  
300& 48.02(3.42)& 21.56(0.98)& 108.77(2.70) \\ 
&600&272.59(5.41)& 61.78(5.56) &323.02(8.12) \\ 
\hline 
\end{tabular}}
\end{table}

\section{Discussion}
In this paper, we aim to provide a general theory for the asymptotic optimality of the classification rules in a large family in the high-dimensional settings with arbitrary number of classes. Our main theorem provides easy-to-check criteria for asymptotic optimality of the classification rules in this family and we establish the corresponding convergence rates as both dimensionality and sample size go to infinity and the number of classes is arbitrary. This general theory is applied to the classic LDA,  the linear programming discriminant rule by \cite{cai2011direct}, and the sparse linear discriminant analysis rule by \cite{Shao-2011}. We extend the latter two methods to the case of multiclass. We establish the asymptotic optimality of the three methods and provide the convergence rates in the high-dimensional settings with arbitrary number of classes. Through simulation study, we demonstrate that the extended methods have good predictive performance when the conditions of these methods are satisfied.

\section*{Acknowledgments}
The second author is supported by NSF DMS 1208786.

\end{document}


\begin{center}
{\LARGE{\bf  Supplementary materials for ``Asymptotic optimality of sparse linear discriminant analysis with arbitrary number of classes''}}
\end{center}
\vskip 2mm

\baselineskip=14pt
\vskip 2mm
\begin{center}
 Ruiyan Luo\\
\vskip 2mm
Division of Epidemiology and Biostatistics, Georgia State University School of Public Health, One Park Place, Atlanta, GA 30303\\
rluo@gsu.edu\\
\hskip 5mm\\
Xin Qi\\
\vskip 2mm
Department of Mathematics and Statistics, Georgia State University, 30 Pryor Street, Atlanta, GA 30303-3083\\
xqi3@gsu.edu \\
\end{center}
  
\vskip 10mm
 \renewcommand{\proof}{\noindent\underline{\bf Proof of Theorem}}

\setcounter{equation}{0}
\setcounter{table}{0}
\baselineskip=18pt

\newpage
\bigskip
 
\section*{Appendix A: Proofs of theorems}
 \renewcommand{\thesubsection}{A.\arabic{subsection}}
 
In this appendix, we provide the proofs of all the theorems in the manuscript.

\subsection{Proof of Theorem  $\ref{theorem_1}$}

 Given a new observation $\mathbf{x}$, $T_{OPT}(\mathbf{x})$ and $T(\mathbf{x})$ denote the classes to which $\mathbf{x}$ is assigned by the rules $T_{OPT}$ and $T$, respectively. We use $P_{\cdot|\mathbf{X}}$ to denote the conditional probability given the  sample $\mathbf{X}$.

\begin{align}
&R_{T}(\mathbf{X})-R_{{OPT}}=(1-R_{{OPT}})-(1-R_{T}(\mathbf{X}))\notag\\
=&\sum_{i=1}^KP\left(T_{OPT}(\mathbf{x})=i|\mathbf{x}\in \text{the $i$th class}\right)P\left(\mathbf{x}\in \text{the $i$th class}\right)\notag\\
&-\sum_{i=1}^KP_{\cdot|\mathbf{X}}\left(T(\mathbf{x})=i|\mathbf{x}\in \text{the $i$th class}\right)P\left(\mathbf{x}\in \text{the $i$th class}\right)\notag\\
=&\frac{1}{K}\sum_{i=1}^K\left[P \left(T_{OPT}(\mathbf{x})=i|\mathbf{x}\in \text{the $i$th class}\right)-P_{\cdot|\mathbf{X}}\left(T(\mathbf{x})=i|\mathbf{x}\in \text{the $i$th class}\right)\right]\notag\\
=&\frac{1}{K}\sum_{i=1}^K\left[\sum_{j\neq i}P_{\cdot|\mathbf{X}}\left(T_{OPT}(\mathbf{x})=i, T(\mathbf{x})=j|\mathbf{x}\in \text{the $i$th class}\right)\label{150}\\
& \left +P_{\cdot|\mathbf{X}}\left(T_{OPT}(\mathbf{x})=i, T(\mathbf{x})=i|\mathbf{x}\in \text{the $i$th class}\right)-P_{\cdot|\mathbf{X}}\left(T(\mathbf{x})=i|\mathbf{x}\in \text{the $i$th class}\right)\right]\notag\\
\le &\frac{1}{K}\sum_{i=1}^K\sum_{j\neq i}P_{\cdot|\mathbf{X}}\left(T_{OPT}(\mathbf{x})=i, T(\mathbf{x})=j|\mathbf{x}\in \text{the $i$th class}\right)\notag\\
= &\frac{1}{K}\sum_{i=1}^K\sum_{j\neq i}P_{\cdot|\mathbf{X}}\left(T(\mathbf{x})=j, T_{OPT}(\mathbf{x})=i|\mathbf{x}\in \text{the $i$th class}\right).\notag
\end{align}  
However in the special case of $K=2$, the right hand side of $\eqref{150}$ is equal to 
\begin{align*}
&\frac{1}{2}\sum_{i=1}^2\left[\sum_{j\neq i}P_{\cdot|\mathbf{X}}\left(T_{OPT}(\mathbf{x})=i, T(\mathbf{x})=j|\mathbf{x}\in \text{the $i$th class}\right)\\
& \left +P_{\cdot|\mathbf{X}}\left(T_{OPT}(\mathbf{x})=i, T(\mathbf{x})=i|\mathbf{x}\in \text{the $i$th class}\right)-P_{\cdot|\mathbf{X}}\left(T(\mathbf{x})=i|\mathbf{x}\in \text{the $i$th class}\right)\right]\\
= &\frac{1}{2}\sum_{i=1}^2\sum_{j\neq i}\left[P_{\cdot|\mathbf{X}}\left(T_{OPT}(\mathbf{x})=i, T(\mathbf{x})=j|\mathbf{x}\in \text{the $i$th class}\right)\\
& \left -P_{\cdot|\mathbf{X}}\left(T_{OPT}(\mathbf{x})=j, T(\mathbf{x})=i|\mathbf{x}\in \text{the $i$th class}\right)\right]\\
= &\frac{1}{2}\sum_{i=1}^2\sum_{j\neq i}\left[P_{\cdot|\mathbf{X}}\left(T(\mathbf{x})=j, T_{OPT}(\mathbf{x})=i|\mathbf{x}\in \text{the $i$th class}\right)\\
& \left -P_{\cdot|\mathbf{X}}\left(T(\mathbf{x})=i, T_{OPT}(\mathbf{x})=j|\mathbf{x}\in \text{the $i$th class}\right)\right]
\end{align*}  
We use $P_i(\cdot)$ to denote the conditional probability $P_{\cdot|\mathbf{X}}\left(\cdot|\mathbf{x}\in \text{the $i$th class}\right)$. Then
\begin{align*}
&P_{\cdot|\mathbf{X}}\left(T(\mathbf{x})=j, T_{OPT}(\mathbf{x})=i|\mathbf{x}\in \text{the $i$th class}\right)=P_i\left(T(\mathbf{x})=j, T_{OPT}(\mathbf{x})=i\right)\\
=& P_i\left(\widehat{\mathbf{a}}_{kj}\trans\boldsymbol{\Sigma}^{-1/2}(\mathbf{x}- \widehat{\mathbf{b}}_{kj})<0, \forall k\neq j, \text{ and } \mathbf{a}_{li}\trans\boldsymbol{\Sigma}^{-1/2}(\mathbf{x}-\mathbf{b}_{li})<0, \forall l\neq i \right)\\
\le & P_i\left(\widehat{\mathbf{a}}_{ij}\trans\boldsymbol{\Sigma}^{-1/2}(\mathbf{x}- \widehat{\mathbf{b}}_{ij})<0, \text{ and } \mathbf{a}_{ji}\trans\boldsymbol{\Sigma}^{-1/2}(\mathbf{x}-\mathbf{b}_{ji})<0 \right)\\
=& P_i\left(\widehat{\mathbf{a}}_{ji}\trans\boldsymbol{\Sigma}^{-1/2}(\mathbf{x}- \widehat{\mathbf{b}}_{ji})>0, \text{ and } \mathbf{a}_{ji}\trans\boldsymbol{\Sigma}^{-1/2}(\mathbf{x}-\mathbf{b}_{ji})<0\right)\\
=& P \left(\widehat{\mathbf{a}}_{ji}\trans\mathbf{Z}>\widehat{\mathbf{a}}_{ji}\trans\boldsymbol{\Sigma}^{-1/2}(\widehat{\mathbf{b}}_{ji}-\boldsymbol{\mu}_i), \text{ and } \mathbf{a}_{ji}\trans\mathbf{Z}<\mathbf{a}_{ji}\trans\boldsymbol{\Sigma}^{-1/2}(\mathbf{b}_{ji}-\boldsymbol{\mu}_i)\bigg\vert\mathbf{X} \right),
\end{align*}  
where $\mathbf{Z}=\boldsymbol{\Sigma}^{-1/2}(\mathbf{x}-\boldsymbol{\mu}_i)\sim N(\mathbf{0},\mathbf{I}_p)$ and independent of the  sample $\mathbf{X}$, and in the fourth line, we use the assumptions $\eqref{24}$. Similarly,
\begin{align*}
&P_{\cdot|\mathbf{X}}\left(T(\mathbf{x})=i, T_{OPT}(\mathbf{x})=j|\mathbf{x}\in \text{the $i$th class}\right)=P_i\left(T(\mathbf{x})=j, T_{OPT}(\mathbf{x})=i\right)\\
=& P \left(\widehat{\mathbf{a}}_{ji}\trans\mathbf{Z}<\widehat{\mathbf{a}}_{ji}\trans\boldsymbol{\Sigma}^{-1/2}(\widehat{\mathbf{b}}_{ji}-\boldsymbol{\mu}_i), \text{ and } \mathbf{a}_{ji}\trans\mathbf{Z}>\mathbf{a}_{ji}\trans\boldsymbol{\Sigma}^{-1/2}(\mathbf{b}_{ji}-\boldsymbol{\mu}_i)\bigg\vert\mathbf{X} \right),
\end{align*}  
 
\newpage

\subsection{Proof of Theorem  $\ref{theorem_7}$}

Given a pair $1\le i\neq j\le K$, we calculate the probability 
$$P \left(\widehat{\mathbf{a}}_{ji}\trans\mathbf{Z}>\widehat{\mathbf{a}}_{ji}\trans\boldsymbol{\Sigma}^{-1/2}(\widehat{\mathbf{b}}_{ji}-\boldsymbol{\mu}_i), \quad\mathbf{a}_{ji}\trans\mathbf{Z}<\mathbf{a}_{ji}\trans\boldsymbol{\Sigma}^{-1/2}(\mathbf{b}_{ji}-\boldsymbol{\mu}_i)\bigg\vert\mathbf{X} \right)$$ 
which is equal to
$P \left(\widehat{\mathbf{a}}_{ji}\trans\mathbf{Z}>\widehat{d}_{ji}, \mathbf{a}_{ji}\trans\mathbf{Z}<d_{ji}\bigg\vert\mathbf{X} \right)$ by the definitions of $d_{ji}$ and $\widehat{d}_{ji}$. Since the subscript $ij$ is fixed during the calculation, for simplicity, we omit it in the following. Moreover, since we calculate the probability with respect to $\mathbf{Z}$, we also omit $ \vert\mathbf{X}$ for simplicity and just regard $\widehat{\mathbf{a}}$ and $\widehat{\mathbf{b}}$ as constants. Note the orthogonal decomposition $\mathbf{a}=t\widehat{\mathbf{a}}+\mathbf{a}_\perp$ and the relationship $d=\frac{1}{2}\|\mathbf{a}\|_2^2$. By the conditions in $\eqref{104}$,
\begin{align}
&\|\mathbf{a}\|^2_2-\|\widehat{\mathbf{a}}\|_2^2=\|\mathbf{a}\|_2^2O_p(s_n),\quad t=1+O_p(s_n), \quad d-\widehat{d}
=\|\widehat{\mathbf{a}}\|_2^2O_p(s_n),\notag\\
\text{we have}\quad &\|\mathbf{a}_\perp\|_2^2=\|\mathbf{a}\|^2_2-t^2\|\widehat{\mathbf{a}}\|_2^2=\|\widehat{\mathbf{a}}\|_2^2O_p(s_n),\quad \widehat{d}
=\|\widehat{\mathbf{a}}\|_2^2(\frac{1}{2}+O_p(s_n)).\label{14}
\end{align} 
We first assume that $\mathbf{a}_\perp\neq \mathbf{0}$. Define
\begin{align*}
& \mathbf{W}=\frac{\widehat{\mathbf{a}}\trans}{\|\widehat{\mathbf{a}}\|_2}\mathbf{Z}\sim N(0,1), \quad \mathbf{V}=-\frac{\mathbf{a}_\perp\trans}{\|\mathbf{a}_\perp\|_2}\mathbf{Z}\sim N(0,1),
\end{align*}
where the distributions are conditional on the  sample $\mathbf{X}$. Since $(\mathbf{W},\mathbf{V})$ is jointly normal and $\widehat{\mathbf{a}}$ and $\mathbf{a}_\perp$ are orthogonal, $\mathbf{W}$ and $\mathbf{V}$ are uncorrelated and hence independent. Let $\phi$ and $\Phi$ are the density and cumulative distribution functions of $N(0,1)$, respectively. Define
\begin{align} 
\eta=\frac{|t\widehat{d}-d|}{t}+\frac{\|\mathbf{a}_\perp\|_2}{t}\sqrt{\log{\left[(\|\mathbf{a}\|^2_2s_n)^{-1}\right]}}.\label{108}
\end{align}
 Then 
\begin{align}
& P\left(\widehat{\mathbf{a}}\trans\mathbf{Z}>\widehat{d}, \mathbf{a}\trans\mathbf{Z}<d\right)=P\left(\widehat{\mathbf{a}}\trans\mathbf{Z}>\widehat{d}, (t\widehat{\mathbf{a}}+\mathbf{a}_\perp)\trans\mathbf{Z}<d\right)\notag\\
=&P\left(\mathbf{W}>\frac{\widehat{d}}{\|\widehat{\mathbf{a}}\|_2},\quad t\|\widehat{\mathbf{a}}\|_2\mathbf{W}-\|\mathbf{a}_\perp\|_2\mathbf{V}<d\right)\notag\\
=&\int_{\frac{\widehat{d}}{\|\widehat{\mathbf{a}}\|_2}}^\infty\phi(w)P\left(\mathbf{V}>\frac{t\|\widehat{\mathbf{a}}\|_2w-d}{\|\mathbf{a}_\perp\|_2}\right)dw=\int_{\frac{\widehat{d}}{\|\widehat{\mathbf{a}}\|_2}}^\infty\phi(w)\left[1-\Phi\left(\frac{t\|\widehat{\mathbf{a}}\|_2w-d}{\|\mathbf{a}_\perp\|_2}\right)\right]dw\notag\\
=&\int_{\frac{\widehat{d}}{\|\widehat{\mathbf{a}}\|_2}}^{\frac{\widehat{d}+\eta}{\|\widehat{\mathbf{a}}\|_2}}+\int_{\frac{\widehat{d}+\eta}{\|\widehat{\mathbf{a}}\|_2}}^\infty \phi(w)\left[1-\Phi\left(\frac{t\|\widehat{\mathbf{a}}\|_2w-d}{\|\mathbf{a}_\perp\|_2}\right)\right]dw\notag\\
\le&\int_{\frac{\widehat{d}}{\|\widehat{\mathbf{a}}\|_2}}^{\frac{\widehat{d}+\eta}{\|\widehat{\mathbf{a}}\|_2}}\phi(w)dw+\int_{\frac{\widehat{d}+\eta}{\|\widehat{\mathbf{a}}\|_2}}^\infty \phi(w)\left[1-\Phi\left(\frac{t\|\widehat{\mathbf{a}}\|_2\frac{\widehat{d}+\eta}{\|\widehat{\mathbf{a}}\|_2}-d}{\|\mathbf{a}_\perp\|_2}\right)\right]dw\notag\\
\le &\frac{\eta}{\|\widehat{\mathbf{a}}\|_2}\phi(\frac{\widehat{d}}{\|\widehat{\mathbf{a}}\|_2})+\left[1-\Phi\left(\frac{t\widehat{d}-d+t\eta}{\|\mathbf{a}_\perp\|_2}\right)\right]\int_{\frac{\widehat{d}}{\|\widehat{\mathbf{a}}\|_2}}^\infty \phi(w)dw.\label{39}
\end{align}  
Since $\widehat{d}/\|\widehat{\mathbf{a}}\|_2=\|\mathbf{a}\|_2(1/2+O_p(s_n))$ and by $\eqref{49}$, $\|\mathbf{a}\|_2$ is bounded below, it follows from the well known inequality:
\begin{align}
(1-\frac{1}{x^2})\phi(x)\le x[1-\Phi(x)]\le \phi(x),\quad \forall x>0,\label{40}
\end{align}  
 that there exists a constant $C_3>0$ independent of $p$ such that with probability converging to 1, 
 \begin{align}
C_3\phi(\frac{\widehat{d}}{\|\widehat{\mathbf{a}}\|_2})\le \frac{\widehat{d}}{\|\widehat{\mathbf{a}}\|_2}[1-\Phi(\frac{\widehat{d}}{\|\widehat{\mathbf{a}}\|_2})].\label{41}
\end{align}  
By $\eqref{14}$, $\eqref{40}$, $\eqref{41}$, and the definition $\eqref{108}$ of $\eta$,  the right hand side of $\eqref{39}$,
\begin{align}
 &\frac{\eta}{\|\widehat{\mathbf{a}}\|_2}\phi(\frac{\widehat{d}}{\|\widehat{\mathbf{a}}\|_2})+\left[1-\Phi\left(\frac{t\widehat{d}-d+t\eta}{\|\mathbf{a}_\perp\|_2}\right)\right][1-\Phi(\frac{\widehat{d}}{\|\widehat{\mathbf{a}}\|_2})]\le\frac{1}{C_3}\frac{\eta}{\|\widehat{\mathbf{a}}\|_2}\frac{\widehat{d}}{\|\widehat{\mathbf{a}}\|_2}[1-\Phi(\frac{\widehat{d}}{\|\widehat{\mathbf{a}}\|_2})]\notag\\
&+\left[1-\Phi\left(\frac{t\widehat{d}-d+|t\widehat{d}-d|+\|\mathbf{a}_\perp\|_2\sqrt{\log{\left[(\|\mathbf{a}\|^2_2s_n)^{-1}\right]}}}{\|\mathbf{a}_\perp\|_2}\right)\right][1-\Phi(\frac{\widehat{d}}{\|\widehat{\mathbf{a}}\|_2})]\notag\\
&\le\frac{1}{C_3}(\frac{1}{2}+O_p(s_n))\eta[1-\Phi(\frac{\widehat{d}}{\|\widehat{\mathbf{a}}\|_2})]+\left[1-\Phi\left(\sqrt{\log{\left[(\|\mathbf{a}\|^2_2s_n)^{-1}\right]}}\right)\right][1-\Phi(\frac{\widehat{d}}{\|\widehat{\mathbf{a}}\|_2})]\notag\\
&\le\frac{1}{C_3}(\frac{1}{2}+O_p(s_n))\eta[1-\Phi(\frac{\widehat{d}}{\|\widehat{\mathbf{a}}\|_2})]+\left[1-\Phi\left(\sqrt{\log{\left[(M_{max}s_n)^{-1}\right]}}\right)\right][1-\Phi(\frac{\widehat{d}}{\|\widehat{\mathbf{a}}\|_2})]\notag\\
&\le \left[\frac{\eta}{C_3}(\frac{1}{2}+O_p(s_n))+\frac{\phi\left(\sqrt{\log{\left[(M_{max}s_n)^{-1}\right]}}\right)}{\sqrt{\log{\left[(M_{max}s_n)^{-1}\right]}}}\right][1-\Phi(\frac{\widehat{d}}{\|\widehat{\mathbf{a}}\|_2})],\label{42}
\end{align}  
where the second inequality from the last is due to $\eqref{49}$. By $\eqref{14}$ and the definition $\eqref{108}$ of $\eta$,
\begin{align}
  &\frac{\eta}{C_3}(\frac{1}{2}+O_p(s_n))+\frac{\phi\left(\sqrt{\log{\left[(M_{max}s_n)^{-1}\right]}}\right)}{\sqrt{\log{\left[(M_{max}s_n)^{-1}\right]}}}\notag \\
	=&\|\mathbf{a}\|_2^2O_p(s_n)+\sqrt{\|\mathbf{a}\|_2^2O_p(s_n)\log{\left[(\|\mathbf{a}\|^2_2s_n)^{-1}\right]}}+O\left(\frac{\exp{\left[-\left(\sqrt{\log{\left[(M_{max}s_n)^{-1}\right]}}\right)^2/2\right]}}{\sqrt{\log{\left[(M_{max}s_n)^{-1}\right]}}} \right)\notag\\
	\le &M_{max}O_p(s_n)+\sqrt{M_{max}O_p(s_n)\log{\left[(M_{max}s_n)^{-1}\right]}}+O\left(\frac{\sqrt{M_{max}s_n}}{\sqrt{\log{\left[(M_{max}s_n)^{-1}\right]}}}\right)\notag\\
	=&O_p\left( \sqrt{M_{max}s_n\log{\left[(M_{max}s_n)^{-1}\right]}} \right)\label{43}.
\end{align}  
Next, we estimate 
\begin{align*}
 \left|[1-\Phi(\frac{\widehat{d}}{\|\widehat{\mathbf{a}}\|_2})]-[1-\Phi(\frac{d}{\|\mathbf{a}\|_2})]\right|=\int_{r_1}^{r_2}\phi(x)dx,
\end{align*}  
where $r_1=\min{(\widehat{d}/\|\widehat{\mathbf{a}}\|_2, d/\|\mathbf{a}\|_2)}$, $r_2=\max{(\widehat{d}/\|\widehat{\mathbf{a}}\|_2, d/\|\mathbf{a}\|_2)}$. By $\eqref{40}$ and  $\eqref{14}$,
\begin{align*}
 &\int_{r_1}^{r_2}\phi(x)dx\le (r_2-r_1)\phi(r_1)=(r_2-r_1)O(r_1[1-\Phi(r_1)])=(r_2-r_1)r_1O([1-\Phi(r_1)])\\
=&\left|\frac{\widehat{d}}{\|\widehat{\mathbf{a}}\|_2}-\frac{d}{\|{\mathbf{a}}\|_2}\right|r_1O([1-\Phi(r_1)])=\left|\|\widehat{\mathbf{a}}\|_2(\frac{1}{2}+O_p(s_n))-
\frac{1}{2}\|{\mathbf{a}}\|_2\right|\frac{1}{2}\|{\mathbf{a}}\|_2O([1-\Phi(r_1)])\\
\le & M_{max}O_p(s_n)O([1-\Phi(r_1)])=O_p(M_{max}s_n)O([1-\Phi(r_1)])
\end{align*}  
Therefore, 
\begin{align}
 &\left|[1-\Phi(\frac{\widehat{d}}{\|\widehat{\mathbf{a}}\|_2})]-[1-\Phi(\frac{d}{\|\mathbf{a}\|_2})]\right|\le O_p(M_{max}s_n)O([1-\Phi(\frac{d}{\|\mathbf{a}\|_2})])=o_p(1)O([1-\Phi(\frac{d}{\|\mathbf{a}\|_2})]),\notag\\
&\text{ and hence } [1-\Phi(\frac{\widehat{d}}{\|\widehat{\mathbf{a}}\|_2})]=[1-\Phi(\frac{d}{\|\mathbf{a}\|_2})](1+o_p(1)),\label{50}
\end{align}  
By $\eqref{39}$, $\eqref{42}$, $\eqref{43}$ and $\eqref{50}$
\begin{align}
& P\left(\widehat{\mathbf{a}}\trans\mathbf{Z}>\widehat{d}, \mathbf{a}\trans\mathbf{Z}<d\right)\le O_p\left( \sqrt{M_{max}s_n\log{\left[(M_{max}s_n)^{-1}\right]}} \right)[1-\Phi(\frac{d}{\|\mathbf{a}\|_2})]\notag\\
&= O_p\left( \sqrt{M_{max}s_n\log{\left[(M_{max}s_n)^{-1}\right]}} \right)P\left( \mathbf{a}\trans\mathbf{Z}>d\right).\label{109}
\end{align}  
Now we consider the case of $\mathbf{a}_\perp=\mathbf{0}$. By similar arguments as those for $\eqref{50}$,
\begin{align}
& P\left(\widehat{\mathbf{a}}\trans\mathbf{Z}>\widehat{d}, \mathbf{a}\trans\mathbf{Z}<d\right)=P\left(\widehat{\mathbf{a}}\trans\mathbf{Z}>\widehat{d}, (t\widehat{\mathbf{a}})\trans\mathbf{Z}<d\right)\label{152}\\
=&P\left(\mathbf{W}>\frac{\widehat{d}}{\|\widehat{\mathbf{a}}\|_2},\quad t\|\widehat{\mathbf{a}}\|_2\mathbf{W}<d\right)\le\left|[1-\Phi(\frac{\widehat{d}}{\|\widehat{\mathbf{a}}\|_2})]-[1-\Phi(\frac{d}{t\|\widehat{\mathbf{a}}\|_2})]\right|\notag\\
 \le& M_{max}O_p(s_n)P\left( \mathbf{a}\trans\mathbf{Z}>d\right)\le O_p\left( \sqrt{M_{max}s_n\log{\left[(M_{max}s_n)^{-1}\right]}} \right)P\left( \mathbf{a}\trans\mathbf{Z}>d\right).\notag
\end{align}  
Combining $\eqref{109}$ and $\eqref{152}$, and note the fact that given the new observation $\mathbf{x}$ belonging to the $i$th class, $\boldsymbol{\Sigma}^{-1/2}(\mathbf{x}-\boldsymbol{\mu}_i)$ have the same distribution as $\mathbf{Z}$, we have
\begin{align}
&P \left(\widehat{\mathbf{a}}_{ji}\trans\mathbf{Z}>\widehat{\mathbf{a}}_{ji}\trans\boldsymbol{\Sigma}^{-1/2}(\widehat{\mathbf{b}}_{ji}-\boldsymbol{\mu}_i),  \mathbf{a}_{ji}\trans\mathbf{Z}<\mathbf{a}_{ji}\trans\boldsymbol{\Sigma}^{-1/2}(\mathbf{b}_{ji}-\boldsymbol{\mu}_i)\right)\label{45}\\
\le &O_p\left( \sqrt{M_{max}s_n\log{\left[(M_{max}s_n)^{-1}\right]}} \right)P \left(\mathbf{a}_{ji}\trans\mathbf{Z}> d_{ji}\right)\notag\\
\le& O_p\left( \sqrt{M_{max}s_n\log{\left[(M_{max}s_n)^{-1}\right]}} \right)\notag\\
&\times P \left(\mathbf{a}_{ji}\trans\boldsymbol{\Sigma}^{-1/2}(\mathbf{x}-\boldsymbol{\mu}_i)> \mathbf{a}_{ji}\trans\boldsymbol{\Sigma}^{-1/2}(\mathbf{b}_{ji}-\boldsymbol{\mu}_i)|\mathbf{x}\in \text{the $i$th class}\right)\notag\\
=&O_p\left( \sqrt{M_{max}s_n\log{\left[(M_{max}s_n)^{-1}\right]}} \right)P\left(\mathbf{a}_{ji}\trans\boldsymbol{\Sigma}^{-1/2}(\mathbf{x}-\mathbf{b}_{ji})>0|\mathbf{x}\in \text{the $i$th class}\right)\notag\\
\le &O_p\left( \sqrt{M_{max}s_n\log{\left[(M_{max}s_n)^{-1}\right]}} \right)P\left(T_{OPT}(\mathbf{x})\notin \text{the $i$th class}|\mathbf{x}\in \text{the $i$th class}\right)\notag\\
\le &O_p\left( \sqrt{M_{max}s_n\log{\left[(M_{max}s_n)^{-1}\right]}} \right)KR_{{OPT}}\notag
\end{align}
 where the $O_p$ term is uniform for all $1\le i\neq j\le K$. Now by Theorem \ref{theorem_1} and $\eqref{45}$,
\begin{align*}
&R_{T}(\mathbf{X})-R_{{OPT}} \\
\le& \frac{1}{K}\sum_{i=1}^K\sum_{j\neq i}\left[P \left(\widehat{\mathbf{a}}_{ji}\trans\mathbf{Z}<\widehat{\mathbf{a}}_{ji}\trans\boldsymbol{\Sigma}^{-1/2}(\widehat{\mathbf{b}}_{ji}-\boldsymbol{\mu}_i),  \quad\mathbf{a}_{ji}\trans\mathbf{Z}> \mathbf{a}_{ji}\trans\boldsymbol{\Sigma}^{-1/2}(\mathbf{b}_{ji}-\boldsymbol{\mu}_i)\right) \right],\notag\\
\le &O_p\left(K^2\sqrt{M_{max}s_n\log{\left[(M_{max}s_n)^{-1}\right]}} \right)R_{{OPT}}\notag
\end{align*}

\newpage

\subsection{Proof of Theorem  $\ref{theorem_8}$}
Similar to the proof of Theorem \ref{theorem_7}, we calculate the probability 
$$P \left(\widehat{\mathbf{a}}_{ji}\trans\mathbf{Z}>\widehat{\mathbf{a}}_{ji}\trans\boldsymbol{\Sigma}^{-1/2}(\widehat{\mathbf{b}}_{ji}-\boldsymbol{\mu}_i),  \mathbf{a}_{ji}\trans\mathbf{Z}<\mathbf{a}_{ji}\trans\boldsymbol{\Sigma}^{-1/2}(\mathbf{b}_{ji}-\boldsymbol{\mu}_i)\bigg\vert\mathbf{X} \right)$$ 
 which is equal to
$P \left(\widehat{\mathbf{a}}_{ji}\trans\mathbf{Z}>\widehat{d}_{ji}, \mathbf{a}_{ji}\trans\mathbf{Z}<d_{ji}\bigg\vert\mathbf{X} \right)$ by the definitions of $d_{ji}$ and $\widehat{d}_{ji}$. As in the proof of Theorem \ref{theorem_7}, we omit the subscript $ij$ and $\vert\mathbf{X} $. Since in the theorem, we assume that $\widehat{\mathbf{a}}$ and $\mathbf{a}$ are not in the same direction, $\mathbf{a}_\perp\neq \mathbf{0}$. $\mathbf{W}$ and $\mathbf{V}$ are the same as in the proof of Theorem \ref{theorem_7}, and let 
\begin{align} 
\delta=-\frac{(t\widehat{d}-d)}{t}+\frac{\|\mathbf{a}_\perp\|_2}{t}.\label{208}
\end{align}
Then by $\eqref{14}$ and $\eqref{186}$
\begin{align}
&\delta\ge  c_5 \sqrt{\|\mathbf{a}\|^2_2s_n}+ o_p\left(\sqrt{\|\mathbf{a}\|^2_2s_n}\right).\label{196}
\end{align}
We calculate
\begin{align}
& P\left(\widehat{\mathbf{a}}\trans\mathbf{Z}>\widehat{d}, \mathbf{a}\trans\mathbf{Z}<d\right)=P\left(\widehat{\mathbf{a}}\trans\mathbf{Z}>\widehat{d}, (t\widehat{\mathbf{a}}+\mathbf{a}_\perp)\trans\mathbf{Z}<d\right)\notag\\
=&P\left(\mathbf{W}>\frac{\widehat{d}}{\|\widehat{\mathbf{a}}\|_2},\quad t\|\widehat{\mathbf{a}}\|_2\mathbf{W}-\|\mathbf{a}_\perp\|_2\mathbf{V}<d\right)\notag\\
=&\int_{\frac{\widehat{d}}{\|\widehat{\mathbf{a}}\|_2}}^\infty\phi(w)P\left(\mathbf{V}>\frac{t\|\widehat{\mathbf{a}}\|_2w-d}{\|\mathbf{a}_\perp\|_2}\right)dw=\int_{\frac{\widehat{d}}{\|\widehat{\mathbf{a}}\|_2}}^\infty\phi(w)\left[1-\Phi\left(\frac{t\|\widehat{\mathbf{a}}\|_2w-d}{\|\mathbf{a}_\perp\|_2}\right)\right]dw\notag\\
\ge &\int_{\frac{\widehat{d}}{\|\widehat{\mathbf{a}}\|_2}}^{\frac{\widehat{d}+\delta}{\|\widehat{\mathbf{a}}\|_2}} \phi(w)\left[1-\Phi\left(\frac{t\|\widehat{\mathbf{a}}\|_2w-d}{\|\mathbf{a}_\perp\|_2}\right)\right]dw\ge\int_{\frac{\widehat{d}}{\|\widehat{\mathbf{a}}\|_2}}^{\frac{\widehat{d}+\delta}{\|\widehat{\mathbf{a}}\|_2}}\phi(w)\left[1-\Phi\left(\frac{t\|\widehat{\mathbf{a}}\|_2\frac{\widehat{d}+\delta}{\|\widehat{\mathbf{a}}\|_2}-d}{\|\mathbf{a}_\perp\|_2}\right)\right]dw\notag\\
\ge&\left[1-\Phi\left(\frac{t\widehat{d}-d+t\delta}{\|\mathbf{a}_\perp\|_2}\right)\right]\int_{\frac{\widehat{d}}{\|\widehat{\mathbf{a}}\|_2}}^{\frac{\widehat{d}+\delta}{\|\widehat{\mathbf{a}}\|_2}} \phi(w)dw\ge [1-\Phi(1)]{\frac{\delta}{\|\widehat{\mathbf{a}}\|_2}}\phi({\frac{\widehat{d}+\delta}{\|\widehat{\mathbf{a}}\|_2}}).\label{189}
\end{align}  
Since
\begin{align*}
 & \phi({\frac{\widehat{d}+\delta}{\|\widehat{\mathbf{a}}\|_2}})/\phi({\frac{\widehat{d}}{\|\widehat{\mathbf{a}}\|_2}})=\exp{(-\frac{2\widehat{d}\delta+\delta^2}{\|\widehat{\mathbf{a}}\|^2_2})}= 1+o_p(1),
\end{align*}  
 and by $\eqref{40}$ and $\eqref{50}$, 
\begin{align}
  &{\frac{\delta}{\|\widehat{\mathbf{a}}\|_2}}\phi({\frac{\widehat{d}}{\|\widehat{\mathbf{a}}\|_2}})\ge {\frac{\delta}{\|\widehat{\mathbf{a}}\|_2}}{\frac{\widehat{d}}{\|\widehat{\mathbf{a}}\|_2}}[1-\Phi(\frac{\widehat{d}}{\|\widehat{\mathbf{a}}\|_2})]\notag\\
	\ge &\left[c_5 \sqrt{\|\mathbf{a}\|^2_2s_n}+ o_p\left(\sqrt{\|\mathbf{a}\|^2_2s_n}\right)\right]\left[\frac{1}{2}+O_p(s_n)\right][1-\Phi(\frac{d}{\|\mathbf{a}\|_2})](1+o_p(1)) \notag\\
	=&\left[\frac{c_5}{2} \sqrt{\|\mathbf{a}\|^2_2s_n}+ o_p\left(\sqrt{\|\mathbf{a}\|^2_2s_n}\right)\right][1-\Phi(\frac{d}{\|\mathbf{a}\|_2})],\label{199}
\end{align}  
we have 
{\allowdisplaybreaks
\begin{align*}
& P \left(\widehat{\mathbf{a}}\trans\mathbf{Z}>\widehat{d}, \mathbf{a}\trans\mathbf{Z}<d\right)\ge \left[\frac{c_5}{2} \sqrt{\|\mathbf{a}\|^2_2s_n}+o_p\left(\sqrt{\|\mathbf{a}\|^2_2s_n}\right)\right][1-\Phi(\frac{d}{\|\mathbf{a}\|_2})]\notag\\
\ge& \left[\frac{c_5}{2} \sqrt{M_{min} s_n}- o_p\left(\sqrt{M_{max}s_n}\right)\right][1-\Phi(\frac{d}{\|\mathbf{a}\|_2})]\notag\\
\ge&\left[\frac{c_5}{2\sqrt{c_4}} \sqrt{M_{max} s_n}- o_p\left(\sqrt{M_{max}s_n}\right)\right] P \left( \mathbf{a}\trans\mathbf{Z}>d\right), 
\end{align*}  
}
where the second line follows from $\eqref{49}$, the last equality is due to $\eqref{186}$ and the $o_p$ term is uniform for all $1\le i\neq j\le K$. Hence,  with probability converging to 1,
 \begin{align}
& P\left(\widehat{\mathbf{a}}\trans\mathbf{Z}>\widehat{d}, \mathbf{a}\trans\mathbf{Z}<d\right)\ge  \frac{c_5}{4\sqrt{c_4}} \sqrt{M_{max} s_n} P\left( \mathbf{a}\trans\mathbf{Z}>d\right).\label{209}
\end{align}  
Now  
\begin{align*}
&\frac{1}{K}\sum_{i=1}^K\sum_{j\neq i}\left[P \left(\widehat{\mathbf{a}}_{ji}\trans\mathbf{Z}>\widehat{\mathbf{a}}_{ji}\trans\boldsymbol{\Sigma}^{-1/2}(\widehat{\mathbf{b}}_{ji}-\boldsymbol{\mu}_i),  \mathbf{a}_{ji}\trans\mathbf{Z}<\mathbf{a}_{ji}\trans\boldsymbol{\Sigma}^{-1/2}(\mathbf{b}_{ji}-\boldsymbol{\mu}_i)\right)\right]\\
\ge &\frac{1}{K}\sum_{i=1}^K\sum_{j\neq i}\left[\frac{c_5}{4\sqrt{c_4}} \sqrt{M_{max} s_n}P\left(\mathbf{a}_{ji}\trans\boldsymbol{\Sigma}^{-1/2}(\mathbf{x}-\mathbf{b}_{ji})>0|\mathbf{x}\in \text{the $i$th class}\right)\right]\notag\\
\ge &\frac{c_5}{4\sqrt{c_4}} \sqrt{M_{max} s_n}\frac{1}{K}\sum_{i=1}^K\left[P\left(\mathbf{a}_{ji}\trans\boldsymbol{\Sigma}^{-1/2}(\mathbf{x}-\mathbf{b}_{ji})>0,\text{ for some } j\neq i|\mathbf{x}\in \text{the $i$th class}\right)\right]\notag\\
=&\frac{c_5}{4\sqrt{c_4}} \sqrt{M_{max} s_n}R_{OPT}.
\end{align*}
 
  \newpage

\subsection{Proof of Theorem  $\ref{theorem_3}$}
Because

 \begin{align}
\left|\widehat{\sigma}_{kl}-(1-\frac{K}{n})\sigma_{kl}\right|=\left|\frac{1}{n}\sum_{i=1}^K\sum_{j=1}^{n_i}(\mathbf{x}^{k}_{ij}-\bar{\mathbf{x}}^{k}_i)(\mathbf{x}^{l}_{ij}-\bar{\mathbf{x}}^{l}_i)-(1-\frac{K}{n})\sigma_{kl}\right|\label{28}\\
=\left|\frac{1}{n}\sum_{i=1}^K\sum_{j=1}^{n_i}(\mathbf{x}^{k}_{ij}-\boldsymbol{\mu}^{k}_i)(\mathbf{x}^{l}_{ij}-\boldsymbol{\mu}^{l}_i)-\frac{1}{n}\sum_{i=1}^K n_i(\bar{\mathbf{x}}^{k}_i-\boldsymbol{\mu}^{k}_i)(\bar{\mathbf{x}}^{l}_i-\boldsymbol{\mu}^{l}_i)-(1-\frac{K}{n})\sigma_{kl}\right|\notag\\
\le  \frac{1}{n}\left|\sum_{i=1}^K\sum_{j=1}^{n_i}\left[(\mathbf{x}^{k}_{ij}-\boldsymbol{\mu}^{k}_i)(\mathbf{x}^{l}_{ij}-\boldsymbol{\mu}^{l}_i)  -\sigma_{kl}\right]\right|+ \frac{1}{n}\left| \sum_{i=1}^K\left[n_i(\bar{\mathbf{x}}^{k}_i-\boldsymbol{\mu}^{k}_i)(\bar{\mathbf{x}}^{l}_i-\boldsymbol{\mu}^{l}_i)-\sigma_{kl}\right] \right|,\notag
\end{align}

where $\mathbf{x}^{k}_{ij}$ denotes the $k$-th coordinate of $\mathbf{x}_{ij}$. Note that both $\mathbf{x}_{ij}-\boldsymbol{\mu}_i$ and $\sqrt{n_i}(\bar{\mathbf{x}}_i-\boldsymbol{\mu}_i)$ have the distributions $N(\mathbf{0}, \boldsymbol{\Sigma})$. By Lemma A.3. in \citet{bickel2008regularized}, we have 

 \begin{align}
P\left(\left|\sum_{i=1}^K\sum_{j=1}^{n_i}\left[(\mathbf{x}^{k}_{ij}-\boldsymbol{\mu}^{k}_i)(\mathbf{x}^{l}_{ij}-\boldsymbol{\mu}^{l}_i)  -\sigma_{kl}\right]\right|>n\nu_1\right)\le C_1\exp(-C_2n\nu_1^2),\label{55}\\
\text{ and }P\left(\max_{k,l}\left|\sum_{i=1}^K\sum_{j=1}^{n_i}\left[(\mathbf{x}^{k}_{ij}-\boldsymbol{\mu}^{k}_i)(\mathbf{x}^{l}_{ij}-\boldsymbol{\mu}^{l}_i)  -\sigma_{kl}\right]\right|>n\nu_1\right)\le C_1p^2\exp(-C_2n\nu_1^2),\notag
\end{align}
for any $\nu_1$ less than a constant $\delta$, where $C_1$, $C_2$ and $\delta$ are constants only depending the upper bound $c_0$ of the eigenvalues of $\boldsymbol{\Sigma}$. By taking $\nu_1=M\sqrt{\frac{\log{p}}{n}}$ for arbitrary $M$, we can obtain 
\begin{align}
\max_{k,l}\frac{1}{n}\left|\sum_{i=1}^K\sum_{j=1}^{n_i}\left[(\mathbf{x}^{k}_{ij}-\boldsymbol{\mu}^{k}_i)(\mathbf{x}^{l}_{ij}-\boldsymbol{\mu}^{l}_i)  -\sigma_{kl}\right]\right|=O_p\left(\sqrt{\frac{\log{p}}{n}}\right).\label{31}
\end{align}
For the second term in $\eqref{28}$, define $Z_{ik}=\sqrt{n_i}(\bar{\mathbf{x}}^{k}_i-\boldsymbol{\mu}^{k}_i)$, for any $1\le i\le K$ and $1\le k\le p$. Then
\begin{align}
 \sum_{i=1}^K\left[n_i(\bar{\mathbf{x}}^{k}_i-\boldsymbol{\mu}^{k}_i)(\bar{\mathbf{x}}^{l}_i-\boldsymbol{\mu}^{l}_i)-\sigma_{kl}\right]=\sum_{i=1}^K\left[Z_{ik}Z_{il}-\sigma_{kl}\right]\notag\\
=\frac{1}{4}\sum_{i=1}^K[(Z_{ik}+Z_{il})^2-(\sigma_{kk}+\sigma_{ll}+2\sigma_{kl})]\notag\\
-\frac{1}{4}\sum_{i=1}^K[(Z_{ik}-Z_{il})^2-(\sigma_{kk}+\sigma_{ll}-2\sigma_{kl})].\label{29}
\end{align}
We will derive the upper bound for the first sum in $\eqref{29}$.  Let 
$$Y_i=(Z_{ik}+Z_{il})^2/(\sigma_{kk}+\sigma_{ll}+2\sigma_{kl})-1.$$

 Then $Y_1,\cdots,Y_K$, are i.i.d. random variables with the distribution $\chi^2_1-1$. We will apply the Bernstein's inequality (see Lemma 2.2.11 in \citet{van1996weak}) to $Y_1+\cdots+Y_K$. We first verify the moment condition required by the Bernstein's inequality. For any positive integer $m\ge 3$,
\begin{align*}
&E[|Y_i|^m]=E[|\chi^2_1-1|^m]\le 2^{m-1}\left(E[|\chi^2_1|^m]+1\right)\\
\le &2^{m}E[|\chi^2_1|^m]=2^m[1\cdot 3\cdot 5\cdots (2m-1)]\le 2^m[2^mm!]\\
=&4^mm!Var(Y_i)/2\le D^{m-2}m!Var(Y_i)/2,
\end{align*}
where $D=64$. It is easy to see that when $m=2$, the inequality $E[|Y_i|^m]\le D^{m-2}m!Var(Y_i)/2$ is also true. Hence,  the moment condition for the Bernstein's inequality holds for $Y_i$. Now by the Bernstein's inequality,
\begin{align*}
&P\left( \left|\sum_{i=1}^K[(Z_{ik}+Z_{il})^2-(\sigma_{kk}+\sigma_{ll}+2\sigma_{kl})]\right|>n\nu_2\right)\\
=&P\left( \left|Y_1+\cdots+Y_K\right|>\frac{n\nu_2}{\sigma_{kk}+\sigma_{ll}+2\sigma_{kl}}\right)\\
\le& 2\exp{(-\frac{1}{2}\frac{x^2}{KVar(Y_1)+Dx})}=2\exp{(-\frac{1}{2}\frac{x^2}{2K+Dx})},
 \end{align*}
where $\nu_2$ can be any positive number and $x=n\nu_2/(\sigma_{kk}+\sigma_{ll}+2\sigma_{kl})$. By the facts that $-x^2/(2K+Dx)$ is a decreasing function for $x>0$, and $\sigma_{kk}+\sigma_{ll}+2\sigma_{kl}$ is less that the two times the largest eigenvalue of $\boldsymbol{\Sigma}$ and hence less than $2c_0$, we have 
\begin{align*}
  \exp{(-\frac{1}{2}\frac{x^2}{2K+Dx})}\le \exp{(-\frac{1}{2}\frac{(n\nu_2)^2}{8c_0^2K+2c_0Dn\nu_2})}
 \end{align*}
\begin{align*}
\text{Therefore, }&P\left( \max_{k,l}\left|\sum_{i=1}^K[(Z_{ik}+Z_{il})^2-(\sigma_{kk}+\sigma_{ll}+2\sigma_{kl})]\right|>n\nu_2\right)\\
&\le 2p^2\exp{\left(-\frac{1}{2}\frac{(n\nu_2)^2}{8c_0^2K+2c_0Dn\nu_2}\right)},
 \end{align*}
Take $\nu_2=M\sqrt{\frac{\log(p)}{n}}$ for arbitrary $M$,  we can obtain

 \begin{align}
2p^2\exp{\left(-\frac{1}{2}\frac{(n\nu_2)^2}{8c_0^2K+2c_0Dn\nu_2}\right)}=2p^2\exp{\left(-\frac{1}{2}\frac{M^2n\log(p)}{8c_0^2K+2c_0DM\sqrt{n\log(p)}}\right)}\notag\\
\le 2p^2\exp{\left(-\frac{1}{2}\frac{M^2n\log(p)}{8c_0^2n+2c_0DM\sqrt{nn}}\right)}=2p^2p^{ -\frac{1}{2}\frac{M^2}{8c_0^2+2c_0DM}},\label{30}
 \end{align}
where we use the facts that $K\le n$ and $\log(p)<n$. As $M$ is large enough, the right hand side of $\eqref{30}$ will be small enough for any $p>1$. Therefore, we have 
\begin{align*}
&\max_{k,l}\frac{1}{n}\left|\sum_{i=1}^K[(Z_{ik}+Z_{il})^2-(\sigma_{kk}+\sigma_{ll}+2\sigma_{kl})]\right|=O_p\left(\sqrt{\frac{\log{p}}{n}}\right).
 \end{align*}
Similarly, we can obtain
\begin{align*}
&\max_{k,l}\frac{1}{n}\left|\sum_{i=1}^K[(Z_{ik}-Z_{il})^2-(\sigma_{kk}+\sigma_{ll}-2\sigma_{kl})]\right|=O_p\left(\sqrt{\frac{\log{p}}{n}}\right).
 \end{align*}
Then by $\eqref{29}$,
\begin{align}
&\max_{k,l}\frac{1}{n}\left| \sum_{i=1}^K\left[n_i(\bar{\mathbf{x}}^{k}_i-\boldsymbol{\mu}^{k}_i)(\bar{\mathbf{x}}^{l}_i-\boldsymbol{\mu}^{l}_i)-\sigma_{kl}\right] \right|=O_p\left(\sqrt{\frac{\log{p}}{n}}\right).\label{32}
 \end{align}
It follows from $\eqref{31}$,  $\eqref{32}$ and $\eqref{28}$,
\begin{align*}
&\max_{k,l} \left|\widehat{\sigma}_{kl}-(1-\frac{K}{n})\sigma_{kl}\right|=O_p\left(\sqrt{\frac{\log{p}}{n}}\right).
\end{align*}
Then we have 
\begin{align*}
&\|\widehat{\boldsymbol{\Sigma}}-(1-\frac{K}{n})\boldsymbol{\Sigma}\|\le \sqrt{\sum_{k,l} \left[\widehat{\sigma}_{kl}-(1-\frac{K}{n})\sigma_{kl}\right]^2}\le \sqrt{p^2\max_{k,l} \left[\widehat{\sigma}_{kl}-(1-\frac{K}{n})\sigma_{kl}\right]^2}\\
\le& p\max_{k,l} \left|\widehat{\sigma}_{kl}-(1-\frac{K}{n})\sigma_{kl}\right| =O_p\left(p\sqrt{\frac{\log{p}}{n}}\right).
\end{align*}
 
\newpage

\subsection{Proof of Theorem  $\ref{theorem_2}$}
To apply Theorem \ref{theorem_7}, we need to verify the conditions in $\eqref{104}$. For simplicity, we omit the subscript $ij$ in the following calculation. Let
$\boldsymbol{\delta}=\boldsymbol{\mu}_j-\boldsymbol{\mu}_i$ and $\widehat{\boldsymbol{\delta}}=\bar{\mathbf{x}}_j-\bar{\mathbf{x}}_i$. Then  
\begin{align}
&\mathbf{a}=\boldsymbol{\Sigma}^{-1/2}\boldsymbol{\delta},\quad \widehat{\mathbf{a}}=\boldsymbol{\Sigma}^{1/2}\widehat{\boldsymbol{\Sigma}}^{-1}\widehat{\boldsymbol{\delta}},\quad d=\frac{1}{2}\|\mathbf{a}\|^2_2=\frac{1}{2}\boldsymbol{\delta}\trans\boldsymbol{\Sigma}^{-1}\boldsymbol{\delta}\label{96}\\
& \widehat{d}=\widehat{\mathbf{a}}\trans\boldsymbol{\Sigma}^{-1/2}(\frac{\bar{\mathbf{x}}_i+\bar{\mathbf{x}}_j}{2}-\boldsymbol{\mu}_i)=\frac{1}{2}\widehat{\boldsymbol{\delta}}\trans\widehat{\boldsymbol{\Sigma}}^{-1}\widehat{\boldsymbol{\delta}}+\widehat{\boldsymbol{\delta}}\trans\widehat{\boldsymbol{\Sigma}}^{-1}(\bar{\mathbf{x}}_i-\boldsymbol{\mu}_i).\notag
\end{align}
By $\eqref{34}$, $\|\widehat{\boldsymbol{\Sigma}}-\boldsymbol{\Sigma}\|=O_p\left(p\sqrt{\frac{\log{p}}{n}}+\frac{K}{n}\right)$. Because $p\sqrt{\log{p}}/\sqrt{n}\le s_n\to 0$, we have $p=o(\sqrt{n})$. Note that  we have assumed  
\begin{align}
K-1\le p.\label{135}
\end{align}
Then
\begin{align}
\frac{K}{n}\le \frac{p+1}{n}=o(\frac{\sqrt{n}}{n})=o\left(p\sqrt{\frac{\log{p}}{n}}\right).\label{122}
\end{align}
 Therefore, we have
\begin{align}
&\|\widehat{\boldsymbol{\Sigma}}-\boldsymbol{\Sigma}\|=O_p\left(p\sqrt{\frac{\log{p}}{n}}+\frac{K}{n}\right)=O_p\left(p\sqrt{\frac{\log{p}}{n}}\right).\label{116}
\end{align}
By $\eqref{95}$,
\begin{align}
&\mathbf{a}\trans\widehat{\mathbf{a}}-\|\mathbf{a}\|_2^2=\boldsymbol{\delta}\trans \widehat{\boldsymbol{\Sigma}}^{-1}\widehat{\boldsymbol{\delta}}-\boldsymbol{\delta}\trans \boldsymbol{\Sigma}^{-1}\boldsymbol{\delta}=\boldsymbol{\delta}\trans \left(\widehat{\boldsymbol{\Sigma}}^{-1}\widehat{\boldsymbol{\delta}}-\boldsymbol{\Sigma}^{-1}\boldsymbol{\delta}\right)\notag\\
=&\boldsymbol{\delta}\trans \left(\widehat{\boldsymbol{\Sigma}}^{-1}-\boldsymbol{\Sigma}^{-1}\right)\widehat{\boldsymbol{\delta}}+\boldsymbol{\delta}\trans\boldsymbol{\Sigma}^{-1}\left( \widehat{\boldsymbol{\delta}}- \boldsymbol{\delta}\right).\label{124}
\end{align}  
Note that $ \widehat{\boldsymbol{\delta}}- \boldsymbol{\delta}=(\bar{\mathbf{x}}_j-\boldsymbol{\mu}_j)-(\bar{\mathbf{x}}_i-\boldsymbol{\mu}_i)$. By Conditions \ref{condition_1} and \ref{condition_3},
\begin{align*}
&c_0^{-1}E[\|\bar{\mathbf{x}}_j-\boldsymbol{\mu}_j\|^2_2]\le E[(\bar{\mathbf{x}}_j-\boldsymbol{\mu}_j)\trans \boldsymbol{\Sigma}^{-1}(\bar{\mathbf{x}}_j-\boldsymbol{\mu}_j)]=\frac{p}{n_j}\le \frac{p}{\min_{1\le i\le K}n_i}\le c_2(\frac{pK}{n}).
\end{align*}  
Hence, for any $t>0$,
\begin{align*}
&P(\max_{1\le j\le K}\|\bar{\mathbf{x}}_j-\boldsymbol{\mu}_j\|_2>t)\le \sum_{j=1}^KP(\|\bar{\mathbf{x}}_j-\boldsymbol{\mu}_j\|_2>t)\le \sum_{j=1}^K\frac{1}{t^2}E[\|\bar{\mathbf{x}}_j-\boldsymbol{\mu}_j\|_2^2]\le \frac{1}{t^2}Kc_2(\frac{pK}{n}).
\end{align*}  
We obtain
\begin{align}
& \max_{1\le j\le K}\|\bar{\mathbf{x}}_j-\boldsymbol{\mu}_j\|_2=O_p\left(K\sqrt{\frac{p}{n}}\right),\text{ and hence } \max_{1\le i\neq j\le K}\|\widehat{\boldsymbol{\delta}}- \boldsymbol{\delta}\|_2=O_p\left(K\sqrt{\frac{p}{n}}\right).\label{123}
\end{align}
By $\eqref{96}$ and $\eqref{123}$,
\begin{align*}
& \|\boldsymbol{\delta}\|_2=\|\mathbf{a}\|_2O_p(1), \quad \|\widehat{\boldsymbol{\delta}}\|_2=\|\mathbf{a}\|_2O_p(1),
\end{align*}
where $O_p(1)$ is uniform for all $1\le i\neq j\le K$. Now by $\eqref{116}$, $\eqref{124}$, $\eqref{123}$ and $\eqref{49}$,
\begin{align}
&\mathbf{a}\trans\widehat{\mathbf{a}}-\|\mathbf{a}\|_2^2=\boldsymbol{\delta}\trans \left(\widehat{\boldsymbol{\Sigma}}^{-1}-\boldsymbol{\Sigma}^{-1}\right)\widehat{\boldsymbol{\delta}}+\boldsymbol{\delta}\trans\boldsymbol{\Sigma}^{-1}\left( \widehat{\boldsymbol{\delta}}- \boldsymbol{\delta}\right)\\
=&\|\mathbf{a}\|_2^2O_p\left(p\sqrt{\frac{\log{p}}{n}}\right)+\|\mathbf{a}\|_2O_p\left(K\sqrt{\frac{p}{n}}\right)=\|\mathbf{a}\|_2^2O_p\left(p\sqrt{\frac{\log{p}}{n}}+\frac{K}{\|\mathbf{a}\|_2}\sqrt{\frac{p}{n}}\right)\notag\\
\le & \|\mathbf{a}\|_2^2O_p\left(p\sqrt{\frac{\log{p}}{n}}+\frac{K}{\sqrt{M_{min}}}\sqrt{\frac{p}{n}}\right)= \|\mathbf{a}\|_2^2O_p(s_n)
\end{align}  
 similarly, $\mathbf{a}\trans\widehat{\mathbf{a}}-\|\widehat{\mathbf{a}}\|_2^2=\|\mathbf{a}\|_2^2O_p(s_n)$. Therefore,
\begin{align}
\|\mathbf{a}\|^2_2-\|\widehat{\mathbf{a}}\|_2^2=\|\mathbf{a}\|_2^2O_p(s_n),  \quad t=\frac{\mathbf{a}\trans\widehat{\mathbf{a}}}{\|\mathbf{a}\|_2^2}=1+O_p(s_n).
\end{align*}  
Moreover, by $\eqref{96}$ and $\eqref{123}$
\begin{align*}
&d-\widehat{d}=\frac{1}{2}\boldsymbol{\delta}\trans \boldsymbol{\Sigma}^{-1}\boldsymbol{\delta}-\frac{1}{2}\widehat{\boldsymbol{\delta}}\trans\widehat{\boldsymbol{\Sigma}}^{-1}\widehat{\boldsymbol{\delta}}-\widehat{\boldsymbol{\delta}}\trans\widehat{\boldsymbol{\Sigma}}^{-1}(\bar{\mathbf{x}}_i-\boldsymbol{\mu}_i)=\|\mathbf{a}\|_2^2O_p(s_n)=\|\widehat{\mathbf{a}}\|_2^2O_p(s_n).
\end{align*}  
Hence we have verified all the conditions in $\eqref{104}$. The theorem follows from Theorem \ref{theorem_7}.
  
\newpage

\subsection{Proof of Theorem  $\ref{theorem_4}$}
We use the idea of the proof of Theorem 1 in \citet{bickel2008covariance}. Let $T_{t_n}$ be the thresholding operator defined in (3) of \citet{bickel2008covariance}. Then 
\begin{align*}
\widetilde{\boldsymbol{\Sigma}}=T_{t_n}\left((1-K/n)^{-1}\widehat{\boldsymbol{\Sigma}}\right),
\end{align*}
and 
\begin{align*}
&\|\widetilde{\boldsymbol{\Sigma}}-\boldsymbol{\Sigma}\|=\|  T_{t_n}\left((1-K/n)^{-1}\widehat{\boldsymbol{\Sigma}}\right)-\boldsymbol{\Sigma}\|\\
\le &\| T_{t_n}\left(\boldsymbol{\Sigma}\right)-\boldsymbol{\Sigma}\|+\| T_{t_n}\left((1-K/n)^{-1}\widehat{\boldsymbol{\Sigma}}\right)-T_{t_n}\left(\boldsymbol{\Sigma}\right)\|.
\end{align*}
For the first term above, by (5.1)
\begin{align}
&\| T_{t_n}\left(\boldsymbol{\Sigma}\right)-\boldsymbol{\Sigma}\|\le \max_{k}\sum_{l=1}^p|\sigma_{kl}|I_{(|\sigma_{kl}|<t_n)}\label{52}\\
=&\max_{k}\sum_{l=1}^p|\sigma_{kl}|^h|\sigma_{kl}|^{1-h}I_{(|\sigma_{kl}|<t_n)}\le \max_{k}\sum_{l=1}^p|\sigma_{kl}|^ht_n^{1-h}I_{(|\sigma_{kl}|<t_n)}\notag\\
\le& t_n^{1-h}\max_{k}\sum_{l=1}^p|\sigma_{kl}|^h =M_1^{1-h}\left(\frac{\log{p}}{n}\right)^{(1-h)/2}C_{h,p}.\notag
\end{align}
On the other hand,
\begin{align}
\left\| T_{t_n}\left((1-K/n)^{-1}\widehat{\boldsymbol{\Sigma}}\right)-T_{t_n}\left(\boldsymbol{\Sigma}\right)\right\|\label{60}\\
\le \max_{k}\sum_{l=1}^p(1-K/n)^{-1}|\widehat{\sigma}_{kl}|I_{[(1-K/n)^{-1}|\widehat{\sigma}_{kl}|\ge t_n, |\sigma_{kl}|<t_n]}\notag\\
+\max_{k}\sum_{l=1}^p |\sigma_{kl}|I_{[(1-K/n)^{-1}|\widehat{\sigma}_{kl}|<t_n, |\sigma_{kl}|\ge t_n]}\notag\\
+\max_{k}\sum_{l=1}^p|(1-K/n)^{-1}\widehat{\sigma}_{kl}-\sigma_{kl}|I_{[(1-K/n)^{-1}|\widehat{\sigma}_{kl}|\ge t_n, |\sigma_{kl}|\ge t_n]}=I+II+III.\notag
\end{align}
By Theorem 4.1 and (5.1),
 \begin{align}
III&\le \max_{k,l}|(1-K/n)^{-1}\widehat{\sigma}_{kl}-\sigma_{kl}|\max_{k}\sum_{l=1}^p |\sigma_{kl}|^h/t_n^h \label{61}\\
&=O_p\left(\left(\frac{\log{p}}{n}\right)^{(1-h)/2}C_{h,p}\right)\notag
\end{align}
For term $II$,

 \begin{align}
II\le\max_{k}\sum_{l=1}^p \left(|(1-K/n)^{-1}\widehat{\sigma}_{kl}-\sigma_{kl}|+|(1-K/n)^{-1}\widehat{\sigma}_{kl}|\right)I_{[(1-K/n)^{-1}|\widehat{\sigma}_{kl}|<t_n, |\sigma_{kl}|\ge t_n]}\notag\\
\le \max_{k,l}|(1-K/n)^{-1}\widehat{\sigma}_{kl}-\sigma_{kl}|\sum_{l=1}^p I_{[|\sigma_{kl}|\ge t_n]}+t_n\max_{k}\sum_{l=1}^pI_{[|\sigma_{kl}|\ge t_n]}\label{62}\\
=O_p\left(\left(\frac{\log{p}}{n}\right)^{(1-h)/2}C_{h,p}\right).\notag
\end{align}
The term 
\begin{align}
 I\le \max_{k}\sum_{l=1}^p |(1-K/n)^{-1}\widehat{\sigma}_{kl}-\sigma_{kl}|I_{[(1-K/n)^{-1}|\widehat{\sigma}_{kl}|\ge t_n, |\sigma_{kl}|< t_n]}\notag\\
+\max_{k}\sum_{l=1}^p |\sigma_{kl}|I_{[|\sigma_{kl}|< t_n]}=IV+V,\label{63}
\end{align}
where by $\eqref{52}$, 
\begin{align}
&V\le M_1^{1-h}\left(\frac{\log{p}}{n}\right)^{(1-h)/2}C_{h,p}.\label{53}
\end{align}
Take a constant $0<\gamma<1$. Then
 \begin{align*}
 &IV=\max_{k}\sum_{l=1}^p |(1-K/n)^{-1}\widehat{\sigma}_{kl}-\sigma_{kl}|I_{[(1-K/n)^{-1}|\widehat{\sigma}_{kl}|\ge t_n, |\sigma_{kl}|< t_n]}\\
&= \max_{k}\sum_{l=1}^p |(1-K/n)^{-1}\widehat{\sigma}_{kl}-\sigma_{kl}|I_{[(1-K/n)^{-1}|\widehat{\sigma}_{kl}|\ge t_n, |\sigma_{kl}|< \gamma t_n]}\\
&\quad +\max_{k}\sum_{l=1}^p |(1-K/n)^{-1}\widehat{\sigma}_{kl}-\sigma_{kl}|I_{[(1-K/n)^{-1}|\widehat{\sigma}_{kl}|\ge t_n, \gamma t_n\le |\sigma_{kl}|< t_n]}\\
\le&\max_{k,l}|(1-K/n)^{-1}\widehat{\sigma}_{kl}-\sigma_{kl}|\max_{k}\sum_{l=1}^p I_\left(|(1-K/n)^{-1}\widehat{\sigma}_{kl}-\sigma_{kl}|>(1-\gamma)t_n\right)\\
&\quad +\max_{k,l}|(1-K/n)^{-1}\widehat{\sigma}_{kl}-\sigma_{kl}|\max_{k}\sum_{l=1}^p  \left(\frac{\sigma_{kl}}{\gamma t_n}\right)^h\\
=&O_p\left(\sqrt{\frac{\log{p}}{n}}\right)\max_{k}\sum_{l=1}^p I_{\left(|(1-K/n)^{-1}\widehat{\sigma}_{kl}-\sigma_{kl}|>(1-\gamma)t_n\right)}+O_p\left(\left(\frac{\log{p}}{n}\right)^{(1-h)/2}C_{h,p}\right).
\end{align*}
By $\eqref{55}$ and $\eqref{30}$ in the proof of Theorem 4.1,
\begin{align*}
 &P\left(\max_{k}\sum_{l=1}^p I_{\left(|(1-K/n)^{-1}\widehat{\sigma}_{kl}-\sigma_{kl}|>(1-\gamma)t_n\right)}>0\right)\\
\le &P\left(\max_{k,l} |(1-K/n)^{-1}\widehat{\sigma}_{kl}-\sigma_{kl}|>(1-\gamma)t_n \right)\\
&\le C_1p^2p^{ -C_2M_2^2}+4p^2p^{ -\frac{1}{2}\frac{M_2^2}{8c_0^2+2c_0DM_2}},
\end{align*}
where $M_2=(1-\gamma)M_1$. If we take $M_1$ large enough, the above probability goes to zero as $p\to \infty$, and 
\begin{align}
&IV=O_p\left(\left(\frac{\log{p}}{n}\right)^{(1-h)/2}C_{h,p}\right).\label{57}
\end{align}
Combining $\eqref{52}$, $\eqref{60}$,$\eqref{61}$,$\eqref{62}$,$\eqref{63}$,$\eqref{53}$,$\eqref{57}$, we have
\begin{align}
&\|\widetilde{\boldsymbol{\Sigma}}-\boldsymbol{\Sigma}\|=O_p\left(\left(\frac{\log{p}}{n}\right)^{(1-h)/2}C_{h,p}\right).\label{66}
\end{align}
Now because $\|\boldsymbol{\Sigma}\|$ is bounded below by $c_0^{-1}$, by $\eqref{66}$, with probability converging to 1, $\|\widetilde{\boldsymbol{\Sigma}}\|$ is bounded below by $c_0^{-1}/2$. Hence, $\|\boldsymbol{\Sigma}^{-1}\|=O(1)$ and $\|\widetilde{\boldsymbol{\Sigma}}^{-1}\|=O_p(1)$, and 
\begin{align*}
&\|\widetilde{\boldsymbol{\Sigma}}^{-1}-\boldsymbol{\Sigma}^{-1}\|=\|\widetilde{\boldsymbol{\Sigma}}^{-1}(\boldsymbol{\Sigma}-\widetilde{\boldsymbol{\Sigma}})\boldsymbol{\Sigma}^{-1}\|\le\|\widetilde{\boldsymbol{\Sigma}}^{-1}\|\|\boldsymbol{\Sigma}-\widetilde{\boldsymbol{\Sigma}}\|\|\boldsymbol{\Sigma}^{-1}\| \notag\\
=&O_p\left(\left(\frac{\log{p}}{n}\right)^{(1-h)/2}C_{h,p}\right).\label{66}
\end{align*}
 
\newpage

\subsection{Proof of Theorem  $\ref{theorem_5}$}
To apply Theorem \ref{theorem_7}, we need to verify the conditions in $\eqref{104}$. For simplicity, we omit the subscript $ij$ in the following calculation. In this particular case, the $\widehat{\mathbf{a}}_{ji}$ and $\widehat{\mathbf{b}}_{ji}$ in Theorem \ref{theorem_1} have the following expressions,
\begin{align}
\widehat{\mathbf{a}}_{ji}= {\boldsymbol{\Sigma}}^{1/2}\widetilde{\boldsymbol{\Sigma}}^{-1}\widetilde{\boldsymbol{\delta}}_{ji}, \quad \widehat{\mathbf{b}}_{ji}=\widetilde{\mathbf{b}}_{ji}=\bar{\mathbf{x}}_1+\frac{\widetilde{\boldsymbol{\delta}}_{j1}+\widetilde{\boldsymbol{\delta}}_{i1}}{2},\quad \forall 1\le i\neq j\le K.\label{71}
\end{align} 
Given $1<j\le K$, as in the proof of Theorem 3 in \citet{Shao-2011}, we have
\begin{align*}
&\|\widetilde{\boldsymbol{\delta}}_{j1}-\boldsymbol{\delta}_{j1}\|^2=O_p(a_n^{2(1-g)}D_{g,p})+O(q_n/n_j),
\end{align*}
where the first term on the right hand side actually does not depend on $j$ because in the definition $\eqref{68}$ of $D_{g,p}$, we have taken the maximum over all pairs of $1\le i\neq j\le K$, and the second term depends on $j$ only through $n_j$ by the same reason. By Condition \ref{condition_3}, we can obtain  
\begin{align*}
&\max_{1<j\le K}\|\widetilde{\boldsymbol{\delta}}_{j1}-\boldsymbol{\delta}_{j1}\|=O_p(\sqrt{a_n^{2(1-g)}D_{g,p}})+O(\sqrt{\frac{Kq_n}{n}})=O_p\left(\sqrt{k_n}\right),
\end{align*}
where $k_n=\max{\{a_n^{2(1-g)}D_{g,p}, Kq_n/n\}}$. Then by our definition of $\widetilde{\boldsymbol{\delta}}_{ji}$, 
\begin{align}
&\max_{1<i\neq j\le K}\|\widetilde{\boldsymbol{\delta}}_{ji}-\boldsymbol{\delta}_{ji}\|=O_p\left(\sqrt{k_n}\right).\label{75}
\end{align}
As in the proof of Theorem 3 in \citet{Shao-2011}, given $j\neq 1$, 
\begin{align*}
\widetilde{\boldsymbol{\delta}}_{j1}\trans\widetilde{\boldsymbol{\Sigma}}^{-1}(\bar{\mathbf{x}}_1-\boldsymbol{\mu}_1)=O_p\left(\sqrt{q_n/n_1}  +\sqrt{C_{h,p}q_n/n}\right)\sqrt{\widetilde{\boldsymbol{\delta}}_{j1}\trans\widetilde{\boldsymbol{\Sigma}}^{-1} \widetilde{\boldsymbol{\delta}}_{j1}},
\end{align*}
where the $O_p$ term does not depend on $j$. Because 
\begin{align*}
&O_p\left(\sqrt{q_n/n_1} +\sqrt{C_{h,p}q_n/n}\right)\sqrt{\widetilde{\boldsymbol{\delta}}_{j1}\trans\widetilde{\boldsymbol{\Sigma}}^{-1} \widetilde{\boldsymbol{\delta}}_{j1}}\\
=&O_p\left(\sqrt{Kq_n/n} +\sqrt{C_{h,p}q_n/n}\right)\sqrt{{\boldsymbol{\delta}}_{j1}\trans{\boldsymbol{\Sigma}}^{-1} {\boldsymbol{\delta}}_{j1}}O_p(1)\\
=&\|\mathbf{a}_{j1}\|O_p\left(\sqrt{(C_{h,p}+K)q_n/n}\right),
\end{align*}
we have $\widetilde{\boldsymbol{\delta}}_{j1}\trans\widetilde{\boldsymbol{\Sigma}}^{-1}(\bar{\mathbf{x}}_1-\boldsymbol{\mu}_1)=\|\mathbf{a}_{j1}\|O_p\left(\sqrt{(C_{h,p}+K)q_n/n}\right)$, and hence
\begin{align}
\widetilde{\boldsymbol{\delta}}_{ji}\trans\widetilde{\boldsymbol{\Sigma}}^{-1}(\bar{\mathbf{x}}_1-\boldsymbol{\mu}_1)=\|\mathbf{a}_{ji}\|O_p\left(\sqrt{(C_{h,p}+K)q_n/n}\right),\label{129}
\end{align}
where the $O_p$ term does not depend on $ij$. Since the following calculation the subscript $ij$ is fixed, for simplicity, we omit it. Then
\begin{align}
& d=\frac{1}{2}\|\mathbf{a}\|^2_2=\frac{1}{2}\boldsymbol{\delta}\trans\boldsymbol{\Sigma}^{-1}\boldsymbol{\delta},\quad  \widehat{d}=\widehat{\mathbf{a}}\trans\boldsymbol{\Sigma}^{-1/2}(\bar{\mathbf{x}}_1+\frac{\widetilde{\boldsymbol{\delta}}_{j1}+\widetilde{\boldsymbol{\delta}}_{i1}}{2}-\boldsymbol{\mu}_i).\label{80} 
\end{align}  
Note that 
{\small\begin{align*}
&\|\boldsymbol{\delta}\|=\|\boldsymbol{\Sigma}^{1/2}\boldsymbol{\Sigma}^{-1/2}\boldsymbol{\delta}\|\le \|\boldsymbol{\Sigma}^{1/2}\|\|\boldsymbol{\Sigma}^{-1/2}\boldsymbol{\delta}\|=\|\boldsymbol{\Sigma}^{1/2}\|\sqrt{\boldsymbol{\delta}\trans\boldsymbol{\Sigma}^{-1/2}\boldsymbol{\delta}}\le\|\mathbf{a}\|_2O_p(1),
\end{align*}}
and $\|\widetilde{\boldsymbol{\delta}}\|=\|\mathbf{a}\|_2O_p(1)$. By  $\eqref{75}$ and Theorem \ref{theorem_4},
{\allowdisplaybreaks
\small\begin{align*}
&|\mathbf{a}\trans\widehat{\mathbf{a}}-\|\mathbf{a}\|_2^2|=|\boldsymbol{\delta}\trans \widetilde{\boldsymbol{\Sigma}}^{-1}\widetilde{\boldsymbol{\delta}}-\boldsymbol{\delta}\trans \boldsymbol{\Sigma}^{-1}\boldsymbol{\delta}|=|\boldsymbol{\delta}\trans \left(\widetilde{\boldsymbol{\Sigma}}^{-1}\widetilde{\boldsymbol{\delta}}-\boldsymbol{\Sigma}^{-1}\boldsymbol{\delta}\right)|\\
=&|\boldsymbol{\delta}\trans \left(\widetilde{\boldsymbol{\Sigma}}^{-1}-\boldsymbol{\Sigma}^{-1}\right)\widetilde{\boldsymbol{\delta}}|+|\boldsymbol{\delta}\trans\boldsymbol{\Sigma}^{-1}\left( \widetilde{\boldsymbol{\delta}}- \boldsymbol{\delta}\right)|\le O_p(\|\mathbf{a}\|_2^2d_n+\|\mathbf{a}\|_2\sqrt{k_n})\\
=&\|\mathbf{a}\|_2^2O_p( d_n+\sqrt{k_n}/\|\mathbf{a}\|_2)\le \|\mathbf{a}\|_2^2O_p( d_n+\sqrt{k_n}/
\sqrt{M_{min}})\\
=&\|\mathbf{a}\|_2^2O_p( \max\{d_n, \sqrt{a_n^{2(1-g)}D_{g,p} }/
\sqrt{M_{min}},\sqrt{{(C_{h,p}+K)q_n/n}}/
\sqrt{M_{min}}\})\\
=&\|\mathbf{a}\|_2^2O_p(b_n)
\end{align*}  }
 and similarly, 
{\allowdisplaybreaks
\small\begin{align*}
&|\mathbf{a}\trans\widehat{\mathbf{a}}-\|\widehat{\mathbf{a}}\|_2^2|=|\boldsymbol{\delta}\trans \widetilde{\boldsymbol{\Sigma}}^{-1}\widetilde{\boldsymbol{\delta}}-\widetilde{\boldsymbol{\delta}}\trans\widetilde{\boldsymbol{\Sigma}}^{-1} \boldsymbol{\Sigma}\widetilde{\boldsymbol{\Sigma}}^{-1}\widetilde{\boldsymbol{\delta}}|\\
&=|\boldsymbol{\delta}\trans \left(\mathbf{I}-\widetilde{\boldsymbol{\Sigma}}^{-1}\boldsymbol{\Sigma}\right)\widetilde{\boldsymbol{\Sigma}}^{-1}\widetilde{\boldsymbol{\delta}}|+|(\boldsymbol{\delta}-\widetilde{\boldsymbol{\delta}})\trans\widetilde{\boldsymbol{\Sigma}}^{-1} \boldsymbol{\Sigma}\widetilde{\boldsymbol{\Sigma}}^{-1}\widetilde{\boldsymbol{\delta}}|\\
&\le O_p(\|\mathbf{a}\|_2^2d_n+\|\mathbf{a}\|_2k_n)=\|\mathbf{a}\|_2^2O_p( d_n+k_n/\|\mathbf{a}\|_2)=\|\mathbf{a}\|_2^2O_p(b_n),
\end{align*}  }
 Therefore,
\begin{align*}
\|\mathbf{a}\|^2_2-\|\widehat{\mathbf{a}}\|_2^2=\|\mathbf{a}\|_2^2O_p(b_n),  \quad t=\frac{\mathbf{a}\trans\widehat{\mathbf{a}}}{\|\mathbf{a}\|_2^2}=1+O_p(b_n).
\end{align*}  
On the other hand, by $\eqref{80}$, $\eqref{75}$ and $\eqref{129}$,
{\allowdisplaybreaks
\small\begin{align*}
&d-\widehat{d}=\frac{1}{2}\boldsymbol{\delta}\trans \boldsymbol{\Sigma}^{-1}\boldsymbol{\delta}- \widetilde{\boldsymbol{\delta}}\trans\widetilde{\boldsymbol{\Sigma}}^{-1}(\bar{\mathbf{x}}_1+\frac{\widetilde{\boldsymbol{\delta}}_{j1}+\widetilde{\boldsymbol{\delta}}_{i1}}{2}-\boldsymbol{\mu}_i)\\
=&\frac{1}{2}\boldsymbol{\delta}\trans \boldsymbol{\Sigma}^{-1}\boldsymbol{\delta}- \widetilde{\boldsymbol{\delta}}\trans\widetilde{\boldsymbol{\Sigma}}^{-1}(\bar{\mathbf{x}}_1+\frac{\widetilde{\boldsymbol{\delta}}_{j1}-\widetilde{\boldsymbol{\delta}}_{i1}}{2}+\widetilde{\boldsymbol{\delta}}_{i1}-\boldsymbol{\mu}_i)\\
=&\frac{1}{2}\boldsymbol{\delta}\trans \boldsymbol{\Sigma}^{-1}\boldsymbol{\delta}-\frac{1}{2}\widetilde{\boldsymbol{\delta}}\trans\widetilde{\boldsymbol{\Sigma}}^{-1}\widetilde{\boldsymbol{\delta}}-\widetilde{\boldsymbol{\delta}}\trans\widetilde{\boldsymbol{\Sigma}}^{-1}(\bar{\mathbf{x}}_1+\widetilde{\boldsymbol{\delta}}_{i1}-\boldsymbol{\mu}_i)\\
=&\frac{1}{2}\boldsymbol{\delta}\trans \boldsymbol{\Sigma}^{-1}\boldsymbol{\delta}-\frac{1}{2}\widetilde{\boldsymbol{\delta}}\trans\widetilde{\boldsymbol{\Sigma}}^{-1}\widetilde{\boldsymbol{\delta}}-\widetilde{\boldsymbol{\delta}}\trans\widetilde{\boldsymbol{\Sigma}}^{-1}(\bar{\mathbf{x}}_1-\boldsymbol{\mu}_1)-\widetilde{\boldsymbol{\delta}}\trans\widetilde{\boldsymbol{\Sigma}}^{-1}(\widetilde{\boldsymbol{\delta}}_{i1}-(\boldsymbol{\mu}_i-\boldsymbol{\mu}_1))\\
=&\frac{1}{2}\boldsymbol{\delta}\trans \boldsymbol{\Sigma}^{-1}\boldsymbol{\delta}-\frac{1}{2}\widetilde{\boldsymbol{\delta}}\trans\widetilde{\boldsymbol{\Sigma}}^{-1}\widetilde{\boldsymbol{\delta}}-\widetilde{\boldsymbol{\delta}}\trans\widetilde{\boldsymbol{\Sigma}}^{-1}(\bar{\mathbf{x}}_1-\boldsymbol{\mu}_1)-\widetilde{\boldsymbol{\delta}}\trans\widetilde{\boldsymbol{\Sigma}}^{-1}(\widetilde{\boldsymbol{\delta}}_{i1}-\boldsymbol{\delta}_{i1})=\|\mathbf{a}\|_2^2O_p(b_n).
\end{align*}  }
All the equalities in $\eqref{104}$ have been proved and the theorem follows from Theorem \ref{theorem_7}.

\newpage
\subsection{Proof of Theorem  $\ref{theorem_6}$}
 We just need to verify the conditions in $\eqref{104}$. In this particular case, 
\begin{align}
&  \mathbf{a}_{ji}=\boldsymbol{\Sigma}^{-1/2}(\boldsymbol{\mu}_j-\boldsymbol{\mu}_i)=\boldsymbol{\Sigma}^{-1/2}\boldsymbol{\delta}_{ji}=\boldsymbol{\Sigma}^{1/2}\boldsymbol{\beta}_{ji}, \quad \widehat{\mathbf{a}}_{ji}={\boldsymbol{\Sigma}}^{1/2}\widehat{\boldsymbol{\beta}}_{ji},\label{170}\\
&d_{ji}=\frac{1}{2}\boldsymbol{\beta}_{ji}\trans \boldsymbol{\Sigma}\boldsymbol{\beta}_{ji},\quad \widehat{d}_{ji}=\widehat{\boldsymbol{\beta}}_{ji}\trans(\frac{\bar{\mathbf{x}}_j+\bar{\mathbf{x}}_i}{2}-\boldsymbol{\mu}_i) \notag
\end{align}  
The following inequality is  a slight revision of (22) in \citet{cai2011direct} and can be proved by the similar argument there. For any $M>0$, we can find a constant $C_1$ such that 
\begin{align*}
&  \max_{1\le j\le K}\max_{1\le k\le p}P\left(|\bar{\mathbf{x}}^k_j-\boldsymbol{\mu}^k_j|\ge C_1\sqrt{\frac{\sigma_{kk}\log{Kp}}{n_j}}\right)\le2p^{-M}. 
\end{align*} 
By Conditions \ref{condition_1} and \ref{condition_3}, $|\sigma_{kk}|$ is bounded by $c_0$ and $n_j\ge n/(c_2K)$, and by $\eqref{135}$,  $K-1\le p$. Hence,
\begin{align*}
&  P\left(\max_{1\le j\le K}\|\bar{\mathbf{x}}_j-\boldsymbol{\mu}_j\|_{\infty}> C_1\sqrt{\frac{c_0c_2K\log{Kp}}{n}}\right)\le2Kpp^{-M}\le 2ppp^{-M}.
\end{align*}
Since $M$ is arbitrary, we have 
\begin{align}
&  \max_{1\le j\le K}\|\bar{\mathbf{x}}_j-\boldsymbol{\mu}_j\|_{\infty}\le O_p\left(\sqrt{\frac{K\log{Kp}}{n}}\right)=O_p\left(\sqrt{\frac{K\log{p}}{n}}\right),\label{134}
\end{align} 
 and then
\begin{align}
&  \max_{1\le i\neq j\le K}\|\widehat{\boldsymbol{\delta}}_{ji}-\boldsymbol{\delta}_{ji}\|_{\infty}\le O_p\left(\sqrt{\frac{K\log{p}}{n}}\right)
.\label{133}
\end{align}  
Similarly, the following inequality is a slight revision of Lemma 2 in \citet{cai2011direct}, for any $M>0$, we can find a constant $C_2$ such that 
\begin{align*}
&  P\left(\max_{1\le i\neq j\le K}\|\overline{\boldsymbol{\Sigma}}\boldsymbol{\beta}_{ji}-\widehat{\boldsymbol{\delta}}_{ji}\|_{\infty}> C_2\sqrt{\frac{M_{max}K\log{p}}{n}}\right)\le 2p^{-M}, 
\end{align*}
We choose the constant $C$ in the definition of $\lambda_n$ to be equal to $C_2$ corresponding to $M=1$ in the above inequality. Then $P\left(\max_{1\le i\neq j\le K}\|\overline{\boldsymbol{\Sigma}}\boldsymbol{\beta}_{ji}-\widehat{\boldsymbol{\delta}}_{ji}\|_{\infty}> \lambda_n\right)\le 2p^{-1}$, and hence 
\begin{align}
  \max_{1\le i\neq j\le K}\|\overline{\boldsymbol{\Sigma}}\boldsymbol{\beta}_{ji}-\widehat{\boldsymbol{\delta}}_{ji}\|_{\infty}=O_p\left(\sqrt{\frac{M_{max}K\log{p}}{n}}\right).\label{136}
\end{align}
Note that in the event $\{\max_{1\le i\neq j\le K}\|\overline{\boldsymbol{\Sigma}}\boldsymbol{\beta}_{ji}-\widehat{\boldsymbol{\delta}}_{ji}\|_{\infty}< \lambda_n\}$, for any $1<j\le K$, $\boldsymbol{\beta}_{j1}$ is in the feasible region of the optimization problem $\eqref{137}$, hence $\|\widehat{\boldsymbol{\beta}}_{j1}\|_1\le \|{\boldsymbol{\beta}}_{j1}\|_1$. Therefore, we have
\begin{align}
  \|\widehat{\boldsymbol{\beta}}_{j1}\|_1\le \|{\boldsymbol{\beta}}_{j1}\|_1+o_p(1),\label{138}
\end{align}
where the $o_p(1)$ term is uniform for all $1<j\le K$. For any $1\le i\neq j\le K$, by $\eqref{133}$, $\eqref{136}$ and $\eqref{138}$,
\begin{align}
&\left|\boldsymbol{\beta}_{j1}\trans  \boldsymbol{\Sigma} \widehat{\boldsymbol{\beta}}_{i1}-\boldsymbol{\beta}_{j1}\trans  \boldsymbol{\Sigma}  \boldsymbol{\beta}_{i1}\right|\le \left|\boldsymbol{\beta}_{j1}\trans\overline{\boldsymbol{\Sigma}} \widehat{\boldsymbol{\beta}}_{i1}-\boldsymbol{\beta}_{j1}\trans  \boldsymbol{\Sigma} \widehat{\boldsymbol{\beta}}_{i1}\right|+\left|\boldsymbol{\beta}_{j1}\trans\overline{\boldsymbol{\Sigma}} \widehat{\boldsymbol{\beta}}_{i1}-\boldsymbol{\beta}_{j1}\trans  \boldsymbol{\Sigma}  \boldsymbol{\beta}_{i1}\right|\label{140}\\
=& \left|\widehat{\boldsymbol{\beta}}_{i1}\trans\left(\overline{\boldsymbol{\Sigma}} \boldsymbol{\beta}_{j1}-\boldsymbol{\delta}_{j1}\right)\right|+\left|\boldsymbol{\beta}_{j1}\trans\left(\overline{\boldsymbol{\Sigma}} \widehat{\boldsymbol{\beta}}_{i1}-\boldsymbol{\delta}_{i1}\right)\right|\le \left|\widehat{\boldsymbol{\beta}}_{i1}\trans\left(\overline{\boldsymbol{\Sigma}} \boldsymbol{\beta}_{j1}-\widehat{\boldsymbol{\delta}}_{j1}\right)\right|+\left|\widehat{\boldsymbol{\beta}}_{i1}\trans\left(\widehat{\boldsymbol{\delta}}_{j1}-\boldsymbol{\delta}_{j1}\right)\right|\notag\\
&+\left|\boldsymbol{\beta}_{j1}\trans\left(\overline{\boldsymbol{\Sigma}} \widehat{\boldsymbol{\beta}}_{i1}-\widehat{\boldsymbol{\delta}}_{i1}\right)\right|+\left|\boldsymbol{\beta}_{j1}\trans\left(\widehat{\boldsymbol{\delta}}_{i1}-\boldsymbol{\delta}_{i1}\right)\right|\le\|\widehat{\boldsymbol{\beta}}_{i1}\|_1\|\overline{\boldsymbol{\Sigma}} \boldsymbol{\beta}_{j1}-\widehat{\boldsymbol{\delta}}_{j1}\|_{\infty}+\|\widehat{\boldsymbol{\beta}}_{i1}\|_1\|\widehat{\boldsymbol{\delta}}_{j1}-\boldsymbol{\delta}_{j1}\|_{\infty}\notag\\
&+\|\boldsymbol{\beta}_{j1}\|_1\|\overline{\boldsymbol{\Sigma}} \widehat{\boldsymbol{\beta}}_{i1}-\widehat{\boldsymbol{\delta}}_{i1}\|_{\infty}+\|\boldsymbol{\beta}_{j1}\|_1\|\widehat{\boldsymbol{\delta}}_{i1}-\boldsymbol{\delta}_{i1}\|_{\infty}\notag\\
\le& O_p\left(\sqrt{\frac{M_{max}K\log{p}}{n}}\right)\max_{1\le i\neq j\le K}\|\boldsymbol{\Sigma}^{-1}\boldsymbol{\delta}_{ji}\|_1.\notag
\end{align}  
Now by $\eqref{140}$,
\begin{align}
&\left|\mathbf{a}_{ji}\trans\widehat{\mathbf{a}}_{ji}-\|\mathbf{a}_{ji}\|_2^2\right|=\left|\boldsymbol{\beta}_{ji}\trans  \boldsymbol{\Sigma} \widehat{\boldsymbol{\beta}}_{ji}-\boldsymbol{\beta}_{ji}\trans  \boldsymbol{\Sigma}  \boldsymbol{\beta}_{ji}\right|=\left|\boldsymbol{\beta}_{ji}\trans  \boldsymbol{\Sigma} (\widehat{\boldsymbol{\beta}}_{ji}-\boldsymbol{\beta}_{ji})\right|\notag\\
=&\left|(\boldsymbol{\beta}_{j1}-\boldsymbol{\beta}_{j1})\trans  \boldsymbol{\Sigma} [(\widehat{\boldsymbol{\beta}}_{j1}-\boldsymbol{\beta}_{j1})-(\widehat{\boldsymbol{\beta}}_{i1}-\boldsymbol{\beta}_{i1})]\right|\le \left|\boldsymbol{\beta}_{j1}\trans  \boldsymbol{\Sigma} (\widehat{\boldsymbol{\beta}}_{j1}-\boldsymbol{\beta}_{j1})]\right|+\left|\boldsymbol{\beta}_{j1}\trans  \boldsymbol{\Sigma} (\widehat{\boldsymbol{\beta}}_{i1}-\boldsymbol{\beta}_{i1})]\right|\notag\\
&+\left|\boldsymbol{\beta}_{i1}\trans  \boldsymbol{\Sigma} (\widehat{\boldsymbol{\beta}}_{j1}-\boldsymbol{\beta}_{j1})]\right|+\left|\boldsymbol{\beta}_{i1}\trans  \boldsymbol{\Sigma} (\widehat{\boldsymbol{\beta}}_{i1}-\boldsymbol{\beta}_{i1})]\right|\le  O_p\left(\sqrt{\frac{M_{max}K\log{p}}{n}}\right)\max_{1\le i\neq j\le K}\|\boldsymbol{\Sigma}^{-1}\boldsymbol{\delta}_{ji}\|_1\notag\\
= &\frac{\|\mathbf{a}_{ji}\|_2^2}{\|\mathbf{a}_{ji}\|_2^2}O_p\left(\sqrt{\frac{M_{max}K\log{p}}{n}}\right)\max_{1\le i\neq j\le K}\|\boldsymbol{\Sigma}^{-1}\boldsymbol{\delta}_{ji}\|_1\notag\\
\le & \frac{\|\mathbf{a}_{ji}\|_2^2}{M_{min}}O_p\left(\sqrt{\frac{M_{max}K\log{p}}{n}}\right)\max_{1\le i\neq j\le K}\|\boldsymbol{\Sigma}^{-1}\boldsymbol{\delta}_{ji}\|_1\le \|\mathbf{a}_{ji}\|_2^2O_p(r_n)
\label{95}
\end{align}  
Next,
\begin{align*}
&\left|\mathbf{a}_{ji}\trans\widehat{\mathbf{a}}_{ji}-\|\widehat{\mathbf{a}}_{ji}\|_2^2\right|=\left|\boldsymbol{\beta}_{ji}\trans  \boldsymbol{\Sigma} \widehat{\boldsymbol{\beta}}_{ji}-\widehat{\boldsymbol{\beta}}_{ji}\trans  \boldsymbol{\Sigma}  \widehat{\boldsymbol{\beta}}_{ji}\right|=\left|\boldsymbol{\delta}_{ji}\trans\widehat{\boldsymbol{\beta}}_{ji}-\widehat{\boldsymbol{\beta}}_{ji}\trans  \overline{\boldsymbol{\Sigma}}  \widehat{\boldsymbol{\beta}}_{ji}\right|+\left| \widehat{\boldsymbol{\beta}}_{ji}\trans  \overline{\boldsymbol{\Sigma}}  \widehat{\boldsymbol{\beta}}_{ji}-\widehat{\boldsymbol{\beta}}_{ji}\trans\boldsymbol{\Sigma}\widehat{\boldsymbol{\beta}}_{ji}\right|,
\end{align*}  
where by the similar argument as in $\eqref{140}$ and $\eqref{95}$, the first term is $\|\mathbf{a}_{ji}\|_2^2O_p(r_n)$. For the second term, by Theorem \ref{theorem_3}, $\eqref{142}$, and $\eqref{138}$,
\begin{align*}
&\left| \widehat{\boldsymbol{\beta}}_{ji}\trans  \overline{\boldsymbol{\Sigma}}  \widehat{\boldsymbol{\beta}}_{ji}-\widehat{\boldsymbol{\beta}}_{ji}\trans\boldsymbol{\Sigma}\widehat{\boldsymbol{\beta}}_{ji}\right|\le \|\widehat{\boldsymbol{\beta}}_{ji}\|_1^2\|\overline{\boldsymbol{\Sigma}}-\boldsymbol{\Sigma}\|_{\infty}=O_p\left(\sqrt{\frac{\log{p}}{n}}\right)\max_{1\le i\neq j\le K}\|\boldsymbol{\Sigma}^{-1}\boldsymbol{\delta}_{ji}\|_1^2\notag\\
\le & \frac{\|\mathbf{a}_{ji}\|_2^2}{M_{min}}O_p\left(\sqrt{\frac{\log{p}}{n}}\right)\max_{1\le i\neq j\le K}\|\boldsymbol{\Sigma}^{-1}\boldsymbol{\delta}_{ji}\|_1^2\le \|\mathbf{a}_{ji}\|_2^2O_p(r_n)
\end{align*} 
Therefore, $\mathbf{a}_{ji}\trans\widehat{\mathbf{a}}_{ji}-\|\widehat{\mathbf{a}}_{ji}\|_2^2=\|\mathbf{a}_{ji}\|_2^2O_p(r_n)$ and 
\begin{align*}
\|\mathbf{a}_{ji}\|^2_2-\|\widehat{\mathbf{a}}_{ji}\|_2^2=\|\mathbf{a}_{ji}\|_2^2O_p(r_n),  \quad t_{ji}=\frac{\mathbf{a}_{ji}\trans\widehat{\mathbf{a}}_{ji}}{\|\mathbf{a}_{ji}\|_2^2}=1+O_p(r_n).
\end{align*}   
On the other hand,  by $\eqref{134}$, $\eqref{138}$ and $\eqref{95}$,
\begin{align*}
&|d_{ji}-\widehat{d}_{ji}|=\left|\frac{1}{2}\boldsymbol{\beta}_{ji}\trans \boldsymbol{\Sigma}\boldsymbol{\beta}_{ji}- \widehat{\boldsymbol{\beta}}_{ji}\trans(\frac{\bar{\mathbf{x}}_j+\bar{\mathbf{x}}_i}{2}-\boldsymbol{\mu}_i)\right|\notag\\
=&\left|\frac{1}{2}\boldsymbol{\beta}_{ji}\trans \boldsymbol{\Sigma}\boldsymbol{\beta}_{ji}- \widehat{\boldsymbol{\beta}}_{ji}\trans(\frac{\boldsymbol{\mu}_j-\boldsymbol{\mu}_i}{2}+\frac{\bar{\mathbf{x}}_j-\bar{\mathbf{x}}_i}{2}-\frac{\boldsymbol{\mu}_j+\boldsymbol{\mu}_i}{2})\right|\\
\le&\left|\frac{1}{2}\boldsymbol{\beta}_{ji}\trans \boldsymbol{\Sigma}\boldsymbol{\beta}_{ji}-\frac{1}{2}\widehat{\boldsymbol{\beta}}_{ji}\trans\boldsymbol{\delta}_{ji}\right|+\left|\left(\frac{\bar{\mathbf{x}}_j+\bar{\mathbf{x}}_i}{2}-\frac{\boldsymbol{\mu}_j+\boldsymbol{\mu}_i}{2}\right)\trans\widehat{\boldsymbol{\beta}}_{ji}\right|\\
=&\left|\mathbf{a}_{ji}\trans\widehat{\mathbf{a}}_{ji}-\|\mathbf{a}_{ji}\|_2^2\right|+\max_{1\le j\le K}\|\bar{\mathbf{x}}_j-\boldsymbol{\mu}_j\|_{\infty}\|\widehat{\boldsymbol{\beta}}_{ji}\|_1\\
\le& \|\mathbf{a}_{ji}\|_2^2O_p(r_n)+O_p\left(\sqrt{\frac{K\log{p}}{n}}\right)\max_{1\le i\neq j\le K}\|\boldsymbol{\Sigma}^{-1}\boldsymbol{\delta}_{ji}\|_1\le \|\mathbf{a}_{ji}\|_2^2O_p(r_n).
\end{align*}
All the conditions in $\eqref{104}$ have been proved.